\documentclass[a4paper]{amsart}
\usepackage[english]{babel}
\usepackage{placeins}
\usepackage[a4paper,margin=1in]{geometry}
\usepackage[utf8]{inputenc}
\usepackage[T1]{fontenc}
\usepackage{graphicx}
\usepackage{amsmath}
\usepackage{amssymb}
\usepackage{amsthm}
\usepackage{mathtools}
\usepackage{listings}
\usepackage{multirow}
\usepackage{enumerate}
\usepackage[autostyle=true,german=quotes]{csquotes}
\usepackage{setspace}

\usepackage[hidelinks]{hyperref}
\usepackage{color}
\usepackage{float}

\definecolor{mygray}{rgb}{0.5,0.5,0.5}
\definecolor{mymauve}{rgb}{0.58,0,0.82}
\definecolor{mygreen}{RGB}{28,172,0} % color values Red, Green, Blue
\definecolor{mylilas}{RGB}{170,55,241}

\newcommand\Keyword[1]{{\sffamily\bfseries #1}}
\newcommand\Do{\enskip\Keyword{do}}

\newcommand\Input{\noindent\Keyword{Input:}\enskip}

\newcommand\Od{\enskip\Keyword{od}}

\newcommand{\R}{{\mathbb R}}
\newcommand{\N}{{\mathbb N}}

\newcommand{\E}{\mathcal{E}}
\newcommand{\T}{\mathcal{T}}

\theoremstyle{remark}

\newtheorem*{rem}{Remark}

\theoremstyle{definition}

\theoremstyle{plain}
\newtheorem{theorem}{Theorem}[section]
\newtheorem{lemma}[theorem]{Lemma}

\newtheorem{remark}[theorem]{Remark}

\addto\captionsenglish{%
}
\addto\captionsngerman{%
}

\usepackage{enumitem}

\renewcommand{\O}{\mathcal O}
\newcommand{\Vertices}{\mathcal V}

\newcommand{\n}{\nu}

\newcommand{\BiL}{\Delta^2}
\newcommand{\JF}[1]{[#1]_E}
\newcommand{\Dn}[1]{{\Dnb #1}}

\newcommand{\DLn}[1]{\partial_{\n} \Delta #1}

\newcommand{\Dnb}[1]{\partial_{\n} #1}
\newcommand{\Dnn}[1]{\partial_{\n\n}^2 #1}

\newcommand{\LFace}{\LEdge}
\newcommand{\LEdge}{{L^2(E)}}

\newcommand{\GCb}{{\Omega}}
\newcommand{\GCSb}{{\Omega}}
\newcommand{\V}{{V}}

\newcommand{\TT}{\mathbb T}
\newcommand{\Tt}{\hat \T}

\newcommand{\uh}{u_h}

\newcommand{\ul}{u_\ell}
\newcommand{\vh}{v_h}
\newcommand{\uht}{\hat \uh}

\newcommand{\Tl}{\T_{\ell}}

\newcommand{\et}{\hat e}

\newcommand{\gh}{g_h}

\renewcommand{\hat}{\widehat}

\renewcommand{\d}{\ensuremath{\ \mathrm{d}}}

\newcommand{\ENorm}[1]{{\left\vert\kern-0.25ex\left\vert\kern-0.25ex\left\vert #1 
    \right\vert\kern-0.25ex\right\vert\kern-0.25ex\right\vert}}
\newcommand{\trb}{\ENorm}

\newcommand{\lh}{\lambda_h}

\newcommand{\huh}{\hat \uh}

\newcommand{\lht}{\widehat \lambda_h}

\newcommand{\lhn}[1]{\lambda_h(#1)}

\renewcommand{\gh}{G\widehat \uh}

\newcommand{\Et}{\widehat{\E}}
\newcommand{\etat}{\widehat{\eta}}
\newcommand{\etal}{\eta_{\ell}}

\newcommand{\W}[1]{{#1}}
\makeatletter
\title[Rate-optimal higher-order conforming FEM for biharmonic
eigenvalue problems]{Rate-optimal 
higher-order adaptive conforming FEM for  biharmonic
eigenvalue problems on polygonal domains}

\author[C.~Carstensen, B.~Gr\"a\ss le]{Carsten Carstensen \and Benedikt Gr\"a\ss le}
\thanks{This work has been supported by the \emph{Deutsche Forschungsgemeinschaft} (DFG) in the Priority Program 1748 \emph{Reliable simulation techniques in solid mechanics: Development of
non-standard discretization methods, 
mechanical and mathematical analysis} 
under the project CA 151/22-2.
The second author is also supported by the DFG under Germany's Excellence Strategy -- The Berlin Mathematics Research Center MATH+ (EXC-2046/1, project ID: 390685689).}

\address[C.~Carstensen, B.~Gr\"a\ss le]{%
	Department of  Mathematics,
	Humboldt-Universit\"at zu Berlin,
	10117 Berlin, Germany; 
\rm cc@math.hu-berlin.de, graesslb@math.hu-berlin.de}
\makeatother

\begin{document}
\maketitle
\begin{abstract}
The a posteriori error analysis of the classical Argyris finite element methods 
dates back to 1996, while the optimal convergence rates of associated
adaptive finite element schemes are established only very recently in 2021. It took a long time 
to realise the necessity of an 
extension of the classical finite element spaces to make them hierarchical.
This paper establishes the novel  adaptive schemes for the biharmonic eigenvalue
problems and provides a mathematical  proof of optimal convergence 
rates towards a simple
eigenvalue and numerical evidence thereof.
This makes the suggested algorithm highly competitive and clearly justifies the
higher computational and implementational costs compared to low-order 
nonconforming schemes. The numerical experiments provide overwhelming evidence that 
higher polynomial degrees pay off with higher convergence rates and underline that adaptive mesh-refining is mandatory. 
Five  computational benchmarks display accurate reference eigenvalues up to 30 digits.
\end{abstract}

\section{Introduction}%
\label{sec:Introduction}
\noindent The conforming Argyris finite element method (FEM) of polynomial order $p\geq5$ allows for an optimally convergent
adaptive scheme that provides guaranteed 
upper eigenvalue bounds  by the min-max principle.

\subsection{State of the art}
The conforming discretisation of the fourth-order problems
with $C^1$ conforming elements like the Hsieh-Clough-Tocher, Bell, or Argyris 
is textbooks material~\cite{Cia:FiniteElementMethod2002}. 
\W{The lowest-order Argyris element consits of quintic polynomials and was also described by
Bell~\cite{Bel:AnalysisThinPlates1968,Bel:RefinedTriangularPlate1969},
Bosshard~\cite{Bos:NeuesVollvertraglichesEndliches1968}, Visser~\cite{Vis:FiniteElementMethod1968},
Withum~\cite{Wit:BerechnungPlattenNach1966}, and
Argyris, Fried, Scharpf~\cite{AFS:TUBAFamilyPlate1968}.}
The a posteriori error control for the Argyris FEM is known for more than two decades and even included in the first
book \cite{Ver:ReviewPosterioriError1996} on a posteriori error analysis. 
But optimal convergence rates of the associated
adaptive algorithm were established only recently for a hierarchical Argyris FEM
\cite{CH:HierarchicalArgyrisFinite2021} for the plain discretisation.
The higher convergence rates of the  Argyris FEM
are rarely visible even in simple
computational benchmarks with the biharmonic equation and a 
right-hand side in $L^2$ in a polygonal bounded Lipschitz
domain $\Omega$.
In contrast, the non-conforming adaptive Morley
FEM is known to be optimal for at least one decade
\cite{HSX:ConvergenceOptimalityAdaptive2012,CGH:DiscreteHelmholtzDecomposition2014}  with an implementation in only 30 lines of MATLAB.

\subsection{Overview}%
\label{sub:Overview}
This paper considers conforming discretisations of eigenvalue problems with hierarchical subspaces
$V(\T_0)\subset V(\T_1)\subset \dots\subset V$ over a 
nested sequence of triangulations $(\Tl)_\ell$ generated by the adaptive algorithm
AFEM of Figure \ref{fig:AFEM}.
The adaptive algorithm refines towards regions of the mesh with a high contribution to the error estimator $\eta$ and
enables optimal convergence
rates.
The key difficulty in the proof of optimal convergence rates for AFEM in the present setting with a hierarchical and
conforming method is the discrete reliability.
The hierarchical structure of the resulting hierarchical Argyris FEM introduced
in~\cite{CH:HierarchicalArgyrisFinite2021} enables guaranteed optimal convergence rates for the source
problem~\cite{CH:HierarchicalArgyrisFinite2021, Gra:OptimalMultilevelAdaptive2022}.
This paper establishes the novel adaptive schemes for biharmonic 
eigenvalue problems and provides a mathematical proof of optimal convergence 
rates towards a simple eigenvalue in Section~\ref{sec:A posteriori}.
The key lemma provides a uniformly improved $L^2$ error control relative to the
$H^2$ energy norm of two discrete eigenvector approximations on two different
but nested triangulations and associated discrete spaces.
The framework 
of elliptic eigenvalue approximations in \cite[Sub.~10.3]{CFPP:AxiomsAdaptivity2014} is 
rewritten and improved. 
Section~\ref{sec:Argyris FEM for eigenvalue problems} establishes the four axioms \eqref{eqn:A1}--\eqref{eqn:A4} for optimal convergence
rates of the hierarchical Argyris eigenvalue solver.
This asserts the high convergence rates of 
the $p$-th order scheme (for $p\ge 5$) in many practical examples. 

The implementation follows \cite{CH:HierarchicalArgyrisFinite2021,Gra:OptimalMultilevelAdaptive2022}
and allows the adaptive computation of all kinds of eigenvalue problems
related to the bi-Laplace operator. 
The first example is the computation of the first reference eigenvalues in more than 30 digits in
Subsection~\ref{sub:high_precision} for the square
and the L-shaped domain in Subsection~\ref{sub:Experiments on the L-shaped domain} with clamped boundary conditions. 
Isospectral domains are well established for the Laplacian and Subsection~\ref{sub:Experiments on isospectral}
investigates this question for the bi-Laplacian: 
There exist isospectral 
domains for the simply-supported plate but those domains have a different
spectrum under clamped  boundary conditions. A more practical example in 
Subsection~\ref{sub:Experiments on a rectangle with a hole} encounters
the computational engineering of a vibrating plate with clamped and free 
boundary conditions. The final example in 
Subsection~\ref{sub:Special constants on triangles} concerns the computation
of optimal constants in basic estimates of numerical analysis: 
Two Friedrichs inequalities and 
interpolation error estimates illustrate that all kind of boundary conditions
can be included.
Optimal convergence rates are visible for all examples and 
thereby clearly justify the higher computational and implementational costs
compared to low-order nonconforming schemes:
The numerical experiments provide overwhelming evidence
that higher polynomial degrees pay off with higher convergence rates
and underline that adaptive mesh-refining is mandatory. 

\subsection{Biharmonic eigenvalue problem and discretisation}%
\label{sub:Biharmonic eigenvalue problem}
Given a bounded Lipschitz domain $\Omega\subset \R^2$ with polygonal boundary
$\partial\Omega$ and outer unit normal $\nu$,
the vector space of admissible functions
\begin{align*}
V\coloneqq\{v\in H^2(\Omega)\ :\ v|_{\partial\Omega}\equiv 0
\equiv \partial_\nu v|_{\partial\Omega}\}
\end{align*}
of continuous functions with square integrable Lebesgue partial derivatives
up to second order and clamped boundary conditions.
The Hessian $Dv\in L^2(\Omega)^{2\times 2}$
defines the energy scalar product 
\[
a(\bullet,\bullet)\coloneqq \left(D^2 \bullet, D^2
\bullet\right)_{L^2(\Omega)}
\]
and with the $L^2$ scalar product 
\[
b(\bullet,\bullet)\coloneqq
\left(\bullet,\bullet\right)_{L^2(\Omega)}.
\]
The weak formulation of the biharmonic eigenvalue problems 
seeks an eigenpair $(\lambda, u)\in\R_+\times V\setminus\{0\}$ with $b(u,u)=1$ and 
\begin{align}\label{eqn:WP}
	a(u,v) = \lambda b(u,v)\quad\text{ for all }v\in V.
\end{align}
Throughout this paper, let $\T$ denote a regular triangulation of the
domain $\Omega$ into compact (non-degenerate) triangles.
The conforming finite element method (FEM) for the discretisation of \eqref{eqn:WP} requires $C^1$ conforming elements to define the
discrete test space $V(\T)$ and seeks a discrete 
eigenpair $(\lh, \uh)\in\R_+\times V(\T)\setminus\{0\}$ with $b(\uh,\uh)=1$ and
\begin{align}\label{eqn:DWP}
	a(\uh,\vh) = \lh b(\uh,\vh)\quad\text{ for all }\vh\in V(\T).
\end{align}
The discrete eigenvalue problem \eqref{eqn:DWP} provides $N\coloneqq\dim V(\T)$ discrete eigenvalues
\begin{align*}
	0<\lhn1\leq\lhn2\leq\dots\leq\lhn N<\infty\eqqcolon\lhn{N+1}\eqqcolon\lhn{N+2}\eqqcolon\cdots.
\end{align*}

\begin{figure}[]\noindent\hspace{-.2em}\fbox{\parbox{.985\textwidth}{\noindent\textbf{AFEM.}\\[.2em]
\Input Initial triangulation $\T_0$,  
$0<\theta<1$, %
and eigenpair
number $j\in\N$
\begin{description}[leftmargin=\widthof{\Keyword{for}\,}]
\item[\Keyword{for}] $\ell=0,1,2,\dots$ \Do
\begin{description}[leftmargin=1em]
\item[\Keyword{Compute}] the $j$-th discrete eigenpair $(\lambda_\ell,\ul)\in \R_+\times V (\Tl)$  with \eqref{eqn:DWP} on $\T\equiv\Tl$
\item[\Keyword{Compute}] $\etal(T):=\eta(T)$ from \eqref{eqn:eta_def} for all $T\in\T\equiv\T_\ell$
\item[\Keyword{Mark}] subset $\mathcal M_\ell\subset\T_\ell$ with 
   $\theta \etal^2(\T) \leq \etal^2(\mathcal{M}_\ell)\coloneqq \sum^{}_{T\in\mathcal{M}}  \eta^2_\ell(T)$ 
    and minimal cardinality $|\mathcal{M}_\ell |\ge 1$
\item[\Keyword{Compute}] smallest NVB refinement $\T_{\ell+1}$ of $\T_\ell$ with
         $\mathcal{M}_\ell\subseteq\T_\ell\setminus\T_{\ell+1}$\Od
\end{description}
\end{description}
\Keyword{Output:} Sequence of triangulations $\T_\ell$, error estimators $\etal$, and eigenpairs 
$(\lambda_\ell,u_\ell)\in\R_+\times V(\T_\ell)$ 
}}
\caption{Adaptive Argyris finite element scheme AFEM}
	\label{fig:AFEM}
\end{figure}

\subsection{Explicit residual-based a  posteriori error estimator}  The local contribution $\eta(T)$ of a triangle  $T\in\T$ with area $|T|>0$ 
and boundary $\partial T$ is the square root of
\begin{align}\label{eqn:eta_def}
\eta^2( T) = |T|^2\|\lh\uh - \BiL \uh\|_{L^2(T)}^2
+|T|^{1/2}\|\JF{\partial_{\nu\nu}^2  u_h}\|_{L^2(\partial T\cap\Omega)}^2+  |T|^{3/2}  \|\JF{\partial_{\nu}\Delta {u_h}}\|_{L^2(\partial T\cap\Omega)}^2.
\end{align}
The jumps $\JF{\bullet}$ apply only along the parts $\partial T\cap\Omega$ of the skeleton inside the domain across an interior edge $E$ 
and only to the derivatives of the piecewise polynomial discrete solution $u_h$; $\Delta $ is the Laplacian and $\Dn\Delta$ is its normal derivative.  

\subsection{Adaptive algorithm AFEM}
The adaptive algorithm depicted in Figure \ref{fig:AFEM} is 
based on the local refinement indicators \eqref{eqn:eta_def}
in a standard (D\"orfler) marking procedure and a
mesh-refining by newest-vertex bisection (NVB) \cite{Ste:CompletionLocallyRefined2008}. 
The NVB defines a set $\TT$ of uniformly shape-regular admissible triangulations.

The theoretical result of this paper is the optimal convergence of the adaptive algorithm
AFEM for sufficiently fine triangulations with mesh-size $h_{\max}\leq \delta$ in $\TT(\delta)$ and a sufficiently small bulk parameter~$\theta$.
\begin{theorem}[optimal rates]\label{thmoptimalrates}
	If $\lambda=\lambda_j$ is a simple eigenvalue, then there exists 
		$0<\delta,\Theta$, that exclusively depend on $\Omega$
	and the shape-regularity of $\TT$ such that for any $0<s<\infty$
	the output $(\Tl)_\ell$ and $(\eta_\ell)_\ell$ of the adaptive algorithm AFEM%
	with bulk parameter $0<\theta<\Theta$ and initial triangulation $\T_0\in\TT(\delta)$ satisfies
	\begin{align}\label{eqccmainresultinpaper}
		\sup_{\ell\in\N_0}
			(1+|\T_\ell| - |\T_0|)^{s}\etal
			\approx\sup_{N\in\N_0}(1+N)^{s}\min \eta(\TT(N)).
	\end{align}
\end{theorem}
	The proof of Theorem~\ref{thmoptimalrates} concludes in 
	Subsection~\ref{sub:proof_optimal_rates}.
\subsection{Outline} The remaining parts of the paper are organised as follows. Section~2 establishes  a rather abstract framework 
for symmetric elliptic eigenvalue problems with Rayleigh quotient in a Gelfand or evolution triple and provides a key lemma for the discrete
stability  on discrete eigenpairs in a 2-level situation. This enables a 2-level  a posteriori error estimation in Lemma
2.2. Those results are essentially known although the framework of \cite[Sub.~10.3]{CFPP:AxiomsAdaptivity2014}
utilises different error terms. 
Section~\ref{sec:Argyris FEM for eigenvalue problems} applies those preliminary results and derives stability (A1),
reduction (A2), discrete reliability (A3), and quasi-orthogonality (A4) for an extended Argyris finite element scheme of degree $p\ge 5$. 
This implies optimal convergence rates \eqref{eqccmainresultinpaper} as the main theoretical  result of this paper. Numerical examples in 
Section~4 provide striking evidence of the recovered higher convergence rates that lead to a rehabilitation of the Argyris finite elements: The higher
implementation costs clearly pay off. It underlines that adaptivity is mandatory and leads to highly accurate results. Five computational benchmarks
provide reference values for the first 10 eigenvalues in high accuracy.

\subsection{General notation}%
\label{sub:General notation}
Standard notation for Lebesgue and Sobolev spaces and their norms 
applies throughout this paper; 
let $H^s(K)$ abbreviate $H^s(\mathrm{int}(K))$ for a (compact) triangle 
$K$ and any $s\in\R$. 
Let $P_k(M)$ denote the spaces of polynomials up to 
degree $k\in\N_0$ on some triangle or edge $M$
with diameter $h_M\in P_0(M)$.
The associated $L^2$ projection $\Pi_{M,k}:L^2(K)\to P_k(M)$ 
is defined by the $L^2$ orthogonality
$(1-\Pi_{M, k})v\perp P_k(M)$ for all $v\in L^2(M)$.
Let
$$P_k(\T)\coloneqq\{p\in L^\infty(\Omega) : p_{|T}\in P_k(T)\text{ for all } T\in T\}$$
denote the space of piecewise polynomials
with respect to a regular triangulation $\T$ of $\Omega$.

\section{Foundations of a posteriori analysis for conforming EVP}%
\label{sec:A posteriori}
This section verifies the axioms of adaptivity of an abstract 
conforming and hierarchical discretisation of a symmetric  elliptic eigenvalue
problem based on the two scalar products $a(\bullet,\bullet)$ and $b(\bullet,\bullet)$.

\subsection{Abstract elliptic eigenvalue problem}%
\label{sub:Conforming EVP analysis}
This subsection considers an evolution or Gelfand triple $(V, H, V^*)$ 
with  the separable infinite-dimensional Hilbert spaces $(V,a)$ and $(H, b)$ 
and their  induced norms
$
\trb{\bullet}\coloneqq a(\bullet, \bullet)^{1/2}$ and $\| \bullet\| \coloneqq b(\bullet,
\bullet)^{1/2}
$
with the compact and dense embeddings 
\begin{align}\label{eqn:evolution_triple}
	V\overset{c}\hookrightarrow H\equiv H^*\overset c\hookrightarrow V^*.
\end{align}
The abstract eigenvalue problem \eqref{eqn:EVP} seeks an eigenpair 
$(\lambda, u)\in\R_+\times V$ with $\|u\|=1$ and
\begin{align}
a(u,v)=\lambda b(u,v)\quad \text{for all }v\in V.\tag{EVP}\label{eqn:EVP}
\end{align}
The well-established spectral theory for the associated compact operators 
(from the compact and dense embeddings in \eqref{eqn:evolution_triple}) 
provides that the countable eigenvalues to \eqref{eqn:EVP} have the
only accumulation point at infinity and, counting their multiplicities, can be enumerated monotonously 
$$
0<\lambda_1\leq\lambda_2\leq\lambda_3\leq\dots \leq \lim_{j\to\infty}\lambda_j=\infty.
$$
The finite dimension $m(j)\in\mathbb N$ of the subordinated 
eigenspaces $E(\lambda_j)=%
\{ u \in V:a(u,\bullet)=\lambda_j b(u,\bullet)\}$ is equal 
to the multiplicity $m(j)\ge 1$ of $\lambda_j$ in the enumeration, i.e.,
$\lambda_k<\lambda_{k+1}=\dots=\lambda_j=\dots=\lambda_{k+m(j)}<\lambda_{k+m(j)+1}$
for each $j\in \{k+1,\dots,k+m(j)\}$ and $k\in \N_0$ with $\lambda_0\coloneqq0$.
Moreover, the corresponding eigenvectors $\phi_1,\phi_2,\dots \in V$
are pairwise $b$-orthonormal as well as $a$-orthogonal, 
\[
a(\phi_j,\psi)=\lambda_j b(\phi_j,\psi) \quad\text{for all }\psi\in V,\quad
a(\phi_j,\phi_k)=\lambda_j  \delta_{jk} ,\quad  \text{and}\quad  b(\phi_j,\phi_k) =\delta_{jk}
\quad\text{for all } j,k\in\N
\]
with the Kronecker delta $ \delta_{jk} $.  We refer to textbooks
\cite{Bre:FunctionalAnalysisSobolev2011,Rud:FunctionalAnalysis1991,Yos:FunctionalAnalysis1995,Zei:NonlinearFunctionalAnalysis1992}
for the functional analysis and the aforementioned facts.

In the 
biharmonic example at hand for a polygonal bounded Lipschitz domain 
$\Omega\subset\R^ 2$, the inner products read
\[
a(u,v):=\int_\Omega D^2 u: D^ 2 v\, dx
\quad  \text{and}\quad
b(u,v):=\int_\Omega  uv\, dx
\]
for $u,v$ in $V:=H^ 2_0(\Omega)$ and~in $H:=L^2(\Omega)$.
The elliptic regularity applies and provides constants $\sigma_\textrm{reg}>1/2$, 
$\sigma:=\min\{1,\sigma_\textrm{reg}\}$,
and  $C_\textrm{reg}$, which exclusively depend on $\Omega$,  
such that, for any $ j\in\N$, the eigenvector 
$\phi_j\in H^{2+\sigma_\textrm{reg}}(\Omega)\cap V$ 
is smooth and satisfies the  regularity estimate 
\[
\| \phi_j \|_{H^{2+\sigma}(\Omega)} \le C_\textrm{reg} \lambda_j 
\]
in the (possibly) broken Sobolev space $H^{2+\sigma}(\Omega)$ with the Sobolev-Slobodeckij 
norm $\|\bullet\|_{H^{2+\sigma}(\Omega)}$~\cite{Gri:EllipticProblemsNonsmooth1985}. 

\subsection{Discrete eigenvalue problem}\label{sub:ConformingEVP}
Consider the family
$\TT$ of regular triangulations (in the sense of Ciarlet) 
of the polygonal bounded Lipschitz domain
$\Omega\subset\R^2$ %
obtained by newest vertex bisection (NVB) based on an initial triangulation $\T_{\textrm{init}}$
\cite{BDD:AdaptiveFiniteElement2004,Ste:CompletionLocallyRefined2008}.
The shape-regularity of the initial triangulation $\T_{\textrm{init}}$ (that matches the domain 
exactly)  transfers  to any triangulation $\T\in\TT$ uniformly. Given any $0<\delta<1$
(resp. $N\in\N_0$) we  abbreviate 
$\TT(\delta)$ (resp.\ $\TT[N]$) as the subset of all $\T\in\TT$ with maximal mesh-size 
$h_{\max}\leq\delta$ (resp.\ with at most $N+|\T_0|$ triangles; here and  
throughout $|\T|$ is the cardinality of $\T$, that is the number of triangles in $\T$). 

Each triangulation $\T\in\TT$ defines  a conforming finite element space $V(\T)\subset V$ 
of  (finite) dimension $N\in\N$ and the associated discrete eigenvalue \eqref{eqn:DWP} problem seeks  some discrete
eigenpair $(\lambda_{h}, u_h)\in\R_+\times V(\T)$ for the number $j$.
There are precisely $N$ 
pairwise $a$-orthogonal  and $b$-orthonormal discrete eigenvectors 
$\phi_h^{(1)} ,\phi_h^{(2)} ,\dots, \phi_h^{(N)}  \in V(\T)$ and those 
form a basis of $V(\T)$.  
Any $j,k=1,\dots,N$ satisfy
\[
a(\phi_h^{(j)} ,\psi_h)=\lambda_h^{(j)} b(\phi_h^{(j)} ,\psi_h) 
\quad\text{for all }\psi_h\in V(\T),\quad
a(\phi_h^{(j)} , \phi_h^{(k)} )=\lambda_h^{(j)} 
 \delta_{jk} ,\quad  \text{and}\quad  b(\phi_h^{(j)} ,\phi_h^{(k)} ) =\delta_{jk}.
\]
Throughout this paper we assume an ascending enumeration 
$0< \lambda_h^{(1)} \le \lambda_h^{(2)} \le \dots
\le \lambda_h^{(N)} $ of the associated eigenvalues counting multiplicities. %
A corollary of the Raleigh-Ritz $\min$-$\max$ principle \cite[Sec.~8]{BO:EigenvalueProblems1991} is the guaranteed upper bound property
\begin{equation}\label{eqccadd1}
\lambda_{k}\le \lambda_h^{(k)}\quad\text{for all }k=1,\dots, N.
\end{equation}

\subsection{A priori  convergence}\label{sub:ConvergenceconformingEVP}
The first convergence results for elliptic eigenvalue problems in question date back to
George Fix in 1973 and the old book \cite{SF:AnalysisFiniteElement1973}  presents an a priori error analysis 
for the Laplacian; but the fairly general arguments apply to the present situation as well;
other standard references include \cite{BO:EigenvalueProblems1991,Bof:FiniteElementApproximation2010}. In the above notation,
let $\lambda:=\lambda_j$ be a simple eigenvalue of number $j\in\N$,
there exists positive constants $\delta_1$ and $\Lambda_\textrm{qo}$ such 
that the following holds
for any triangulation $\T\in\TT(\delta_1)$ with maximal mesh-size at most $\delta_1$.
The dimension $N$ of $V(\T)$ is at least $j$ and 
we simply drop the index $j$ and abbreviate the discrete eigenpair 
$(\lambda_h^{(j)}, \phi_h^{(j)})$ as $(\lambda_h, u_h)$ and the exact eigenpair $(\lambda_j, \phi_j)$
as $(\lambda, u)$ with a possible change of signs in 
$u_h=\pm  \phi_h^{(j)}$ so that $b(u,u_h)\geq0$ and 
\begin{equation}\label{eqccadd2}
\textrm{LB}:=\trb{ u-u_h}^ 2 + \lambda_h-\lambda +\| u-u_h\|^ 2
\le  \Lambda_\textrm{qo}  \min_{v_h\in V(\T)}\trb{ u-v_h}^ 2=\O(h_{\max}^{2\sigma} ).
\end{equation}
Notice that $ \lambda_h-\lambda\ge 0 $ (by \eqref{eqccadd1}) implies that all terms 
of the  lower bound  \textrm{LB} in \eqref{eqccadd2} are non-negative.

Let us summarise the smallness assumptions on the mesh-sizes required so far.
First some $\delta_1>0$ guarantees that the dimension $N:=\dim V(\T)$
of $V(\T)$ is at least $j\in\N$ so that a discrete eigenpair $(\lambda_h,u_h)$
of number $j\le N$ exists whenever $\T\in\TT(\delta_1)$.
 Second, we shall use that  \textrm{LB} in \eqref{eqccadd2} is small and we 
quantify this (with some possibly smaller $\delta_1$) by 
\begin{equation}\label{eqccadd2a}
 \textrm{LB}< 1/2 \quad\text{for all } \T\in \TT(\delta_1).
\end{equation}
This implies in particular that $b(u,u_h)= 1-1/2\, \| u-u_h\|^2\geq3/4$
(from the binomial theorem and $\|u\|=1=\|u_h\|$) is at least  $3/4$
and the sign of $u_h=\pm  \phi_h^{(j)}$ is uniquely determined by that of $u$. 

We need some further smallness of $\textrm{LB}$ for uniform separation of 
eigenvalues.  
The resolution of the spectral gap 
$\lambda_{j-1}<\lambda\equiv\lambda_j<\lambda_{j+1}$ of the simple eigenvalue 
$\lambda$ with the number $j\in\N$ 
(and the convention  $\lambda_0:=\lambda_h^ {(0)}:=0$)
by the discrete counterparts 
$\lambda_h^ {(j-1)}\le \lambda_h^ {(j)}\equiv\lambda_h\le \lambda_h^ {(j+1)}$ requires
\begin{equation}\label{eqccadd2b}
\lambda_{j-1}\le \lambda_h^ {(j-1)}<\lambda_j\equiv\lambda\le \lambda_h^ {(j)}\equiv\lambda_h<\lambda_{j+1}\le  \lambda_h^ {(j+1)},
\end{equation}
where the positive gaps $\lambda- \lambda_h^ {(j-1)}$ and $\lambda_{j+1}-\lambda_h$
are bounded below by a positive constant. This follows from \eqref{eqccadd2} for sufficiently 
small mesh-sizes  $h_{\max}$ and leads to some positive $ \delta_2\le \delta_1$ such that
\eqref{eqccadd2b} holds for all $\T\in\TT(\delta_2)$.
Consequently, the lower bounds of the weighted gaps provides the boundedness 
\begin{equation}\label{eqccadd2c}
\max\{ 
\frac{ \lambda_j }{\lambda_j-\lambda_h^ {(j-1)} },
\frac{ \lambda_h^ {(j)} }{\lambda_h^ {(j+1)}-\lambda_h^ {(j)} }
\}\le  \Lambda_\textrm{sep}
\end{equation}
with some positive separation constant $ \Lambda_\textrm{sep}<\infty$ that 
is independent of $\T\in\TT(\delta_2)$. Separation conditions frequently 
appear in an eigenvalue convergence analysis and we note that 
\eqref{eqccadd2c} implies in particular that
\[
 \sup_{\T\in\TT(\delta_2)}
 \max_{k\ne j}\frac{ \lambda_h^ {(j)} }{|\lambda_h^ {(k)} -\lambda_h^{(j)}|}\le  \Lambda_\textrm{sep}
 \quad\text{and}\quad
  \sup_{\T\in\TT(\delta_2)}
 \max_{k\ne j}\frac{ \lambda_j }{|\lambda_h^ {(k)} -\lambda_j|}
 \le  \Lambda_\textrm{sep}
\]
(with the abbreviation $ \max_{k\ne j}$ for the maximum 
over all $k\in\{1,\dots, N\}\setminus\{j\}$
with $N:=\dim V(\T)$).
In fact, the monotonicity and $\lambda\le\lambda_h$ 
imply the missing part
$
{ \lambda_h }/(\lambda_h -\lambda_h^{(j-1)})\le { \lambda}/(\lambda -\lambda_h^{(j-1)})
\le  \Lambda_\textrm{sep}
$ of the first estimate;
the proof of the second  is analogous.

There are extensions of \eqref{eqccadd2} 
to multiple eigenvalues with reduced convergence \cite{SF:AnalysisFiniteElement1973},
but then the notation is significantly more technical and the analog analysis of the cluster
algorithm \cite{Gal:OptimalAdaptiveFEM2015} remains for future research.

\subsection{Two-level notation}\label{sub:2levelnotationConvergenceconformingEVP}
Throughout this paper, let
$j\in\N$ denote the fixed number of a simple eigenvalue $\lambda\equiv \lambda_j$ and assume
that the  positive constants
$\delta_2\le \delta_1$, $\Lambda_\textrm{qo}$, and $\Lambda_\textrm{sep}$
satisfy \eqref{eqccadd2}-\eqref{eqccadd2c}.
Consider an admissible refinement $\Tt$
of $\T\in\TT(\delta_2) $ and let $(\lambda_h, u_h)\in \R^ +\times V(\T)$ resp.\
 $(\widehat \lambda_h, \widehat u_h)\in \R^ +\times V(\widehat\T)$ denote 
 the respective discrete eigenpair number $j$ of \eqref{eqn:DWP}. %
 Recall that  
  \eqref{eqccadd2}--\eqref{eqccadd2a} uniquely determines
  (the orientation of) those discrete eigenvectors.
 Another corollary
 of the Rayleigh-Ritz $\min$-$\max$ principle \cite[Sec.~8]{BO:EigenvalueProblems1991} (on the algebraic level) 
 and  $ V(\T)\subset V(\widehat\T)$ reads 
 $\widehat \lambda_h\le \lambda_h$ and so
$
 \lambda\le \widehat \lambda_h\le \lambda_h
$
in addition to \eqref{eqccadd2b}.

The key to any a priori or a posteriori error analysis is the uniform smallness of the $H$ norm compared to the $V$ norm of the
error term $\widehat e:=\widehat u_h-u_h$. 

\begin{lemma}[key]\label{prop:stab}
There exists a  constant $\Lambda_\textrm{stab}$, that exclusively depends on 
$\Lambda_\textrm{sep},\Omega,\sigma, \delta_2$, and $\TT$, but not on 
$\T\in\TT(\delta_2)$ with  maximal mesh-size $h_{\max}\le \delta_2$, 
such that
\begin{align*}
\| \uht-\uh \| \leq\Lambda_\textrm{stab}  h_{\rm{max}}^\sigma \trb{ \uht-\uh}.
\end{align*}
\end{lemma}

\begin{proof}
The lemma is indicated in \cite{CFPP:AxiomsAdaptivity2014}
for the Laplace eigenvalue problem,
but the proof is merely outlined. The  more detailed   proof below has eight steps.
\emph{The first step}  verifies a 2-level separation condition 
\begin{equation}\label{eqccadd5}
 \sup_{\T\in\TT(\delta_2)}
 \max_{k\ne j}\frac{ \widehat \lambda_h }{|\lambda_h^ {(k)} -\widehat \lambda_h|}\le  \Lambda_\textrm{sep}.
\end{equation}
The proof of  \eqref{eqccadd5} follows by monotonicity and the relations 
$\lambda\leq\widehat \lambda_h\le \lambda_h$ in the 
$2$-level notation:
\[
\frac{ \widehat\lambda_h  }{\lambda_h^ {(j+1)} - \widehat\lambda_h }
\le 
\frac{ \lambda_h  }{\lambda_h^ {(j+1)} - \lambda_h }
\le  \Lambda_\textrm{sep}
 \quad\text{and}\quad
 \frac{ \widehat\lambda_h  }{ \widehat\lambda_h -\lambda_h^ {(j-1)}}
\le 
 \frac{ \lambda  }{ \lambda -\lambda_h^ {(j-1)}}
\le  \Lambda_\textrm{sep}.
\]
{\em The second step} concerns  
 the Galerkin projection $G\in L(V)$ onto $ V(\T)$ with
\begin{align}\label{eqn:G_def}
 Gv\in V(\T) \quad\text{and}\quad 	a(v-Gv,\vh)=0\qquad\text{for all }(v,\vh) \in V\times \V(\T).
\end{align}
The approximation properties are well established and a standard duality argument
reveals 
\begin{equation}\label{eqccadd6}
\| v-Gv\|\le \Lambda_\textrm{dual} h_{\max}^ \sigma  \trb{ v-Gv} 
\qquad\text{for all }v \in V
\end{equation}
with a universal constant $\Lambda_\textrm{dual}$ that exclusively depends on 
$\sigma,\Omega,$ and $\TT$, but is independent of $\T\in\TT$ and $v$. Further
details can be found in
\cite{Cia:FiniteElementMethod2002,SF:AnalysisFiniteElement1973,EG:FiniteElementQuasiinterpolation2017,BS:MathematicalTheoryFinite2008}.
\newline\emph{The third step} introduces the $H$ projection 
$G \widehat u_h-\gamma u_h$
of $G \widehat u_h$ onto the span
of $u_h$ with the Fourier coefficient $\gamma:=b(G \widehat u_h, u_h)$.
Since the discrete vector $G \widehat u_h-\gamma u_h\in V(\T)$ belongs to the 
coarse finite element space $ V(\T)$ of dimension $N$, 
it has a Fourier expansion with the 
eigenvectors
$\phi_h^ {(k)}$ and coefficients
$\alpha_k:=b(G \widehat u_h-\gamma u_h, \phi_h^ {(k)})\in\R $ for $k=1,\dots,N$.
The orthogonality of the discrete eigenvectors  implies $\alpha_j=0$ and  
$\alpha_k=b(G \widehat u_h, \phi_h^ {(k)}) $ for  $k\ne j$.
Consequently, the finite  Fourier series 
\[
G \widehat u_h-\gamma u_h=\sum_{k\ne j} \alpha_k \phi_h^ {(k)}
\]
leads to the $H$ norm identity 
\[
\| G \widehat u_h-\gamma u_h\|^ 2= \sum_{k\ne j} \alpha_k^2.
\]
{\em The fourth step} performs the critical spectral analysis with 
 the discrete eigenvalue problem \eqref{eqn:DWP} for an eigenpair 
$(\lh^ {(k)} , \phi_h^ {(k)} )$ of number $k\ne j$ with 
$\alpha_k=b(G \widehat u_h, \phi_h^ {(k)}) $. 
Therefore,  \eqref{eqn:DWP} and  \eqref{eqn:G_def}   provide 
\[
\lambda_h^ {(k)}\alpha_k=a(G \widehat u_h, \phi_h^ {(k)})=a(\widehat u_h, \phi_h^ {(k)}).
\]
Since $(\widehat \lambda_h,\widehat u_h)$ is an eigenpair on the fine level and 
$V(\T)\subset V(\widehat\T)$, 
$a(\widehat u_h, \phi_h^ {(k)})=\widehat \lambda_h b(\widehat u_h, \phi_h^ {(k)})$. 
Thus
\begin{align*}
(\lambda_h^ {(k)}-\widehat\lambda_h)\; \alpha_k=\lht\;		b(\huh-\gh, \phi_h^ {(k)})
=\lht \beta_k
\end{align*}
with the abbreviation $\beta_k:=b(  \widehat u_h- G \widehat u_h,  \phi_h^ {(k)})\in\R $,
 $k\ne j$. This reads $ \alpha_k=\lht \beta_k/(\lambda_h^ {(k)}-\widehat\lambda_h)$ 
 such that
\[
|\alpha_k|\le  \Lambda_\textrm{sep} | \beta_k|\quad\text{for all }k\ne j
\] 
follows from  \eqref{eqccadd5} and establishes the intermediate bound  
\[
\| G \widehat u_h-\gamma u_h\|^ 2= \sum_{k\ne j} \alpha_k^2
\le \Lambda_\textrm{sep}^ 2 \sum_{k\ne j} \beta_k^ 2.
\]
{\em Step five} starts with a Bessel inequality  for the Fourier coefficients
$\beta_k=b(  \widehat u_h- G \widehat u_h,  \phi_h^ {(k)})$, namely
\[
 \sum_{k\ne j} \beta_k^ 2\le \| \widehat u_h- G \widehat u_h\|^ 2.
\]
The synopsis with the previous step and thereafter  \eqref{eqccadd6} verify
\[
\| G \widehat u_h-\gamma u_h\|\le \Lambda_\textrm{sep}\| \widehat u_h- G \widehat u_h\|
\le  \Lambda_\textrm{sep}
 \Lambda_\textrm{dual} h_{\max}^ \sigma  \trb{  \widehat u_h-G \widehat u_h} .
\]
Since the Galerkin projection \eqref{eqn:G_def} is a best-approximation,
$ \trb{  \widehat u_h-G \widehat u_h} \le  \trb{  \widehat u_h-u_h} 
= \trb{  \widehat e}$ results in
\begin{equation}\label{eqccadd7}
\| G \widehat u_h-\gamma u_h\|\le  \Lambda_\textrm{sep} \Lambda_\textrm{dual} h_{\max}^ \sigma  
\trb{  \widehat e}.
\end{equation}
{\em Step six} clarifies $\beta:=b( \widehat u_h,  u_h)\ge 0$. The binomial theorem and
$\|\widehat  u_h\|=1=\|u_h\|$ reveal 
\[
2(1-\beta) =  \| \widehat e\|^ 2 \le \left(\|u-\widehat u_h\|+\|u- u_h\|\right)^2\le 2
\]
with a triangle inequality and  $\|u-\widehat u_h\|, \|u-u_h\|\le 2^ {-1/2}$
from \eqref{eqccadd2a} in the last step. Hence $0\le \beta$.

{\em Step seven} proves $\|\widehat e \| \le \sqrt{2}  \|\huh-\gamma\uh\|$ by elementary 
algebra and $\|\huh\|=1=\|\uh\|$. Indeed,
direct calculations with $0\le \beta\coloneqq b(\huh, \uh)\le 1$ in the last estimate result in
\[
\|\widehat e \|^2/2 = 1-\beta=\|\huh - \beta \uh\|^2/(1+\beta)\le \|\huh - \beta \uh\|^2.
\]
This and
$
 \|\huh - \beta \uh\| = \min_{t\in\R}\|\huh - t\uh\|
$
from $\|\uht-t\uh\|^2=1-2t\beta+t^2$ read
\begin{equation}\label{eqccadd8}
2^{-1/2}  \|\widehat e \| \le  \min_{t\in\R}\|\huh - t\uh\|\le \|\huh - \gamma \uh\|.
\end{equation}
{\em Step eight} concludes the proof. 
Recall \eqref{eqccadd8} and employ a triangle inequality to infer
\[
2^{-1/2}  \|\widehat e\| \le  \|\huh - \gamma \uh\|\le 
 \|\huh -G\huh\|+ \|G\huh - \gamma \uh\|.
\]	
The terms $ \|\huh -G\huh\|$ and $\|G\huh - \gamma \uh\|$
are controlled by \eqref{eqccadd6}--\eqref{eqccadd7}.  The combination with 
$ \trb{ \huh-G\huh} \le  \trb{ \widehat e}$ (from best-approximation of $G$)
leads to
\[
2^{-1/2}  \|\widehat e \| \le 
   (1+\Lambda_\textrm{sep}) \Lambda_\textrm{dual} h_{\max}^ \sigma  
\trb{  \widehat e}.
\]
This is the assertion for  $\Lambda_\textrm{stab}:= \sqrt{2} (1+\Lambda_\textrm{sep}) 
 \Lambda_\textrm{dual}$.
\end{proof}

\subsection{Discrete reliability analysis}
Given the 2-level notation of Subsection~\ref{sub:2levelnotationConvergenceconformingEVP} 
with the fine and coarse discrete eigenpair 
$(\lht, \uht)\in\R_+\times V(\widehat\T)$ and
$(\lh, \uh)\in\R_+\times V(\T)$, respectively,
 and the error $\widehat e:=\uht - \uh\in V(\widehat\T)$,
consider the linear bounded functional $\textrm{Res}\in V^*$ on the coarse level by
\begin{equation}\label{eqccadd9}
\textrm{Res}(\varphi):= \lambda_h b(u_h,\varphi)-a(u_h,\varphi)
\quad\text{for all }\varphi\in V.
\end{equation} 
Observe that the solution property \eqref{eqn:DWP} for the coarse
eigenpair  $(\lambda_h,u_h)\in\R_+\times V(\T)$ implies $\textrm{Res}(v_h)=0$ for all $v_h\in V(\T)$, written
 $ V(\T)\subset \ker \textrm{Res}$. The latter property enables a classical 
 explicit residual-based a posteriori error analysis and allows below 
 the derivation of  computable 
 reliable upper bounds of 
 \[
  \| \textrm{Res}\|_{V(\widehat\T)^*}:=\sup_{\widehat\vh\in V(\widehat\T)\setminus\{0\}} 
 \textrm{Res}(\widehat\vh)/\trb{\widehat\vh}\le
 \trb{ \textrm{Res}}_*:= \| \textrm{Res}\|_{V^*}:=\sup_{v\in V\setminus\{0\}} 
 \textrm{Res}(v)/\trb{v}.
 \] 
Given any parameter $0<\kappa<1$, set  
$\delta_3:=\min\{\delta_2, (\kappa/(\lambda_{j+1} \Lambda^2_\textrm{stab})^{1/(2\sigma)}\}$ and suppose $\T\in\TT(\delta_3)$. %

\begin{lemma}[2-level a posteriori error estimation]\label{lem:err_Res_est}
Under the 2-level notation with  $\T\in\TT(\delta_3)$, it holds
\begin{align*}
(a)\qquad
(1-\kappa) \trb{\et}^2\le\textrm{Res}(\et)\le  \lh-\lht  \le \trb{\et}^2
\quad \text{ and }\quad (b) \quad 
(1-\kappa) \trb{\et}\le   \| \textrm{Res}\|_{V(\widehat\T)^*} \le  2\trb{ \widehat e}.
\end{align*}
\end{lemma}

\begin{proof}[Proof of Lemma \ref{lem:err_Res_est}.a]
The role of the residual becomes clear in a little calculation with the fine
eigenpair $(\widehat\lambda_h,\widehat u_h)$ and the error $\et\in V(\widehat\T)$,
\begin{align}\label{eqn:e_Res_id}
\trb{\et}^2&=a(\uht - \uh, \et) = \lht b(\uht, \et) - a(\uh, \et)
		=\textrm{Res}(\et) + b(\lht\uht  -\lh\uh, \et).
	\end{align}
It remains to control the nonlinear contribution $ b(\lht\uht  -\lh\uh, \et)$ with elementary
algebra and the key lemma.
The binomial theorem and $\|\uht\|=1=\|\uh\|$ imply $b(\uht+\uh,
\uht-\uh)=0$. This and
$2(\lht\uht  -\lh\uh) = (\lht+\lh)\et + (\lht-\lh)(\uht+\uh)$ verify
\begin{align}\label{eqccadd9a}
	2b(\lht\uht  -\lh\uh, \et) = b((\lht+\lh)\et+ (\lht -\lh)(\uht+ \uh),\et)=(\lht +\lh)\|\et\|^2.
\end{align}
Therefore, \eqref{eqn:e_Res_id} reads
\begin{equation}\label{eqccadd10}
\trb{\et}^2=\textrm{Res}(\et) + (\lht +\lh)/2\; \|\et\|^2.
\end{equation}
The additional assumption $\T\in\TT(\delta_3)$ provides 
$h_{\max}^{2\sigma}\lambda_{j+1} \Lambda^2_\textrm{stab}\le\kappa$, 
so that $ (\lht +\lh)/2\le \lambda_{j+1}$ by \eqref{eqccadd2b} leads
in  the key Lemma~\ref{prop:stab} to 
\begin{equation}\label{eqccadd11}
(\lht +\lh)/2 \| \widehat e \|^2 \leq\kappa   \trb{ \widehat e}^2.
\end{equation}
The combination of  \eqref{eqccadd10}--\eqref{eqccadd11} 
results in 
$ (1-\kappa) \trb{\et}^2\le \textrm{Res}(\et) \le  \trb{\et}^2$
as part of the assertions.
The remaining statements follow with the
well known  Pythagorean theorem for eigenpairs
\begin{align}\label{eqccadd12}
\trb{\et}^2= \lh-\lht+\lht\|\et\|^2.
\end{align}
This 2-level version
 is stated and proven in \cite[Lem.\ 10.6]{CFPP:AxiomsAdaptivity2014} for the 
Laplace operator, but the direct algebraic proof for the situation at hand is the 
same and hence not repeated.
Recall  the key Lemma~\ref{prop:stab} for 
$\lambda_{j+1} \| \widehat e \|^2 \leq\kappa   \trb{ \widehat e}^2$ and 
 $ \lht \le \lambda_{j+1}$ (from \eqref{eqccadd2b}) to deduce with \eqref{eqccadd12} that
 \[
 0\le \trb{\et}^2- (\lh-\lht)= \lht\|\et\|^2\le \kappa \trb{\et}^2.
 \]
This is another part of the assertions, namely
$ (1-\kappa) \trb{\et}^2\le \lh-\lht\le \trb{\et}^2 $.
The combination of 
  \eqref{eqccadd10} and \eqref{eqccadd12}  results in
 \(
 \textrm{Res}(\et) =(\lh-\lht)(1- \|\et\|^2/2)  \le \lh-\lht
 \)
 and concludes the proof of $(a)$.
 \end{proof}
 
 \begin{proof}[Proof of Lemma \ref{lem:err_Res_est}.b]
 {\em Step one analyses the Riesz representation $\widehat v$  of $ \textrm{Res}$.}
 Let $\widehat v\in V(\widehat\T)$ denote the Riesz representation of 
 $ \textrm{Res}=a(\widehat v,\bullet)$ in the Hilbert space $( V(\widehat\T),a)$ and 
 deduce $G\widehat v=0$
 for the Galerkin projection \eqref{eqn:G_def} from
 \[
	 \trb{ G\W{\widehat v} }^2=a(\W{\widehat v},G\W{\widehat v})=\textrm{Res}(G\W{\widehat v})=0 
\quad\text{by }G\W{\widehat v}\in V(\T)\subset \ker\textrm{Res}.
 \]
The duality estimate  \eqref{eqccadd6}  and 
$h_{\max}^{\sigma}\le\Lambda_\textrm{stab}^{-1} \sqrt{\kappa/\lambda_{j+1} }$ from 
$\T\in\TT(\delta_3)$ in the last step provide
\[
\|\widehat v\|=\| \widehat v-G\widehat  v\|\le 
\Lambda_\textrm{dual} h_{\max}^ \sigma  \trb{ \widehat  v} 
\le
\Lambda_\textrm{dual}\Lambda_\textrm{stab}^{-1} \sqrt{\kappa/\lambda_{j+1} } \trb{ \widehat  v} .
\] 
 Recall $\Lambda_\textrm{stab}= \sqrt{2} (1+\Lambda_\textrm{sep}) 
 \Lambda_\textrm{dual}$ from the key Lemma~\ref{prop:stab} to recast this as
 \begin{align}\label{eqccadd12b}
 \|\widehat v\|
 \le (1+\Lambda_\textrm{sep})^ {-1}
\sqrt{\kappa/(2\lambda_{j+1}) } 
  \trb{ \widehat  v} .
 \end{align}
 {\em Step two controls the nonlinear term $ \| \lht\uht  -\lh\uh\|$.} 
The arguments, which led above  to \eqref{eqccadd9a},
also verify %
\begin{align*}
	4 \| \lht\uht  -\lh\uh\|^ 2  - (\lht+\lh)^ 2 \|\et\|^2 =(\lh -\lht)^2\| \uht+ \uh\|^ 2
	\le 4 (\lh -\lht)^2.
\end{align*}
Recall \eqref{eqccadd11} %
to rewrite this with part (a) as 
\begin{align*}
	 \| \lht\uht  -\lh\uh\|^ 2  \le \left(  (\lht+\lh)\kappa/2 + \lh -\lht)\right)   \trb{ \widehat e}^2  
	 \le \lambda_{j+1} (1+\kappa/2)  \trb{ \widehat e}^2 
\end{align*}
with  $\lh\le \lambda_{j+1}$ (from \eqref{eqccadd2b}) in the last step.
In other words,
\begin{align}\label{eqccadd13}
	 \| \lht\uht  -\lh\uh\|  \le \sqrt{ \lambda_{j+1} (1+\kappa/2) } \trb{ \widehat e}.
\end{align}
 {\em Step three concludes the proof.} 
The Riesz  representation $\widehat v$ with isometry  
$ \trb{ \textrm{Res}}_{V(\widehat\T)^*}=\trb{\widehat v}$  satisfies
 \[
 \trb{\widehat v}^2=\textrm{Res}(\widehat v)
 =b(\lh\uh,\widehat v)-a(\uh,\widehat v)
 =b(\lh\uh-\widehat \lh\widehat \uh,\widehat v)+a(\widehat e,\widehat v)
 \]
 with the fine
 eigenpair $(\widehat\lambda_h,\widehat u_h)$ and the error $\et\in V(\widehat\T)$
 in the last step. %
 Cauchy inequalities %
 and the estimates  \eqref{eqccadd12b}--\eqref{eqccadd13} reveal
 \[
  \trb{\widehat v}^2\le \left( 1+  (1+\Lambda_\textrm{sep})^ {-1}
 \sqrt{ (1+\kappa/2)\kappa/2 }
  \right)   \trb{ \widehat e} \trb{ \widehat  v} .
 \]
 This proves a sharper version of the asserted estimate
  $ \trb{\widehat v}\le  2\trb{ \widehat e}$. The remaining assertion follows from  part (a)
  and the definition of the operator norm  $ \| \textrm{Res}\|_{V(\widehat\T)^*}$: %
 \(
  (1-\kappa) \trb{ \widehat e}^2\le \textrm{Res}(\widehat e)
  \le    \| \textrm{Res}\|_{V(\widehat\T)^*}  \trb{ \widehat e} .
\)
 \end{proof}

The emphasis of the a posteriori error analysis in this  subsection is 
on the universal constants for a large class of triangulations, i.e., to allow
$\delta_3$ as large as possible. The other extreme when $h_{\rm max}\leq\delta_3$ tends to
zero (e.g., from uniform refinement) remains undisplayed, but is worth a comment. 

\begin{rem}[asymptotic exactness]
If the maximal mesh-size $h_{\max}$ of $\T$ becomes smaller and tends to zero,  the quantities $  \textrm{Res}(\et)$, $ \lh-\lht$,  $\trb{\et}^2$, 
and  $\| \textrm{Res}\|_{V(\widehat\T)^*}$ are asymptotically exact in 
the sense that their 
pairwise quotients converge towards one. This follows essentially from the analysis in this subsection with the key
 Lemma~\ref{prop:stab} from $\kappa\to 0$ as $\delta_3\to 0$ and has been observed before in Laplace eigenvalue problems
 \cite{CG:OscillationfreeAdaptiveFEM2011}.  
\end{rem}

\section{Argyris FEM for eigenvalue problems}%
\label{sec:Argyris FEM for eigenvalue problems}
This section introduces the notation for the standard and hierarchical Argyris FEM.

\subsection{Triangulations and Argyris finite element spaces}%
\label{sub:Triangulation}
The Argyris finite element of degree $p=5,6,7,\dots$ is well established and could be defined on a regular 
triangulation $\T$ with set of vertices $\Vertices$ of the polygonal bounded Lipschitz domain $\Omega\subset\R^2$ by
\begin{align}\label{eqn:A_std}
	A_p(\T)\coloneqq\left\{v_h\in P_p(\T)\cap V  : D^2v_h\text{ is continuous at every } z\in\Vertices\right\}.
\end{align}
This Argyris finite element space from \cite{AFS:TUBAFamilyPlate1968,BS:InteriorPenaltyMethods2005} 
is nowadays textbook topic
\cite{Bra:FiniteElementsTheory2007,BS:MathematicalTheoryFinite2008,BC:FiniteElementMethods2017,Cia:FiniteElementMethod2002}
and we apply the hierarchical 
 Argyris finite element method from \cite{CH:HierarchicalArgyrisFinite2021}. Given a sequence of successive 
 one-level refinements  $\T_0,\T_1,\T_2,\dots$ of the initial triangulation $\T_0$
by the newest-vertex bisection (NVB) 
 (where each triangle is divided in at most four sub-triangles via bisection of refinement edges), the spaces 
 \begin{align}\label{eqn:A_ext}
	V(\T_\ell) \coloneqq A_p(\T_0) + A_p(\T_1) + \dots + A_p(\T_\ell)\quad\text{for any }\ell\in\N_0
\end{align}
are obviously conforming and nested. It can be shown that   $V(\T_\ell) $ exclusively depends on $\T_\ell$ (and does not depend on the 
the sequence   $\T_0,\T_1,\T_2,\dots,\T_\ell$). This and a characterisation of the nodal basis in the extended Argyris finite element space  
$V(\T_\ell)$ are provided in
\cite{CH:HierarchicalArgyrisFinite2021,Gra:OptimalMultilevelAdaptive2022,CG:ManuscriptPreparation2024}. 

Let $\TT$ denote the set of all triangulations $\T_\ell$ generated in the NVB and given 
 $\T\in \TT$  some refinement  $\widehat\T\in \TT(\T)$ thereof, adopt the 2-level notation of 
 Subsection~\ref{sub:2levelnotationConvergenceconformingEVP} with discrete eigenpairs 
 $(\lambda_h, u_h)\in \R^ +\times V(\T)$ resp. $(\widehat \lambda_h, \widehat u_h)\in \R^ +\times V(\widehat\T)$ 
of fixed number $j$ of \eqref{eqn:DWP}; set $\widehat e:= \widehat u_h- u_h\in  V(\widehat\T)$  and 
\[
\delta(\T,\widehat\T):=\trb{\widehat e} .
\]
The subsequent subsections provide some properties \eqref{eqn:A1}--\eqref{eqn:A4} sufficient for optimal convergence rates of AFEM
\cite{CFPP:AxiomsAdaptivity2014,CR:AxiomsAdaptivitySeparate2017,CH:HierarchicalArgyrisFinite2021}. 

\subsection{Proof of stability (A1)} The error estimator \eqref{eqn:eta_def} 
is based on the coarse level of the 2-level notation of 
Subsection~\ref{sub:2levelnotationConvergenceconformingEVP}. Its counterpart on the
fine level with the discrete eigenpair
 $(\widehat \lambda_h, \widehat u_h)\in \R^ +\times V(\widehat\T)$ 
is  $\widehat\eta^2(\T)\coloneqq \sum_{T\in\widehat\T} \widehat\eta^2( T)$
based on the triangulation $\widehat\T\in\TT(\delta_2)$ and the contribution 
\begin{align}\label{eqn:eta_deffine}
\widehat\eta^2( T) = |T|^2\|\widehat\lh\widehat\uh - \BiL \widehat\uh\|_{L^2(T)}^2
+|T|^{1/2}\|\JF{\partial_{\nu\nu}^2  \widehat u_h}\|_{L^2(\partial T\cap\Omega)}^2
+  |T|^{3/2}\|\JF{\Dn\Delta  {\widehat u_h}} \|_{L^2(\partial T\cap\Omega)}^2
\end{align}
for each $T\in \widehat\T$. Recall that $\widehat\T$ is a refinement of $\T$, whence
$\T\cap\widehat\T$ is the (possibly empty) set of coarse and fine triangles. 
The sum convention defines $\eta(\T\cap\widehat\T)$
(resp. $\widehat\eta(\T\cap\widehat\T)$) as the square root of
 $\eta^2(\T\cap\widehat\T)\coloneqq \sum_{T\in\T\cap\widehat\T} \eta^2( T)$
(resp.   $\widehat\eta^2(\T\cap\widehat\T)\coloneqq 
\sum_{T\in\T\cap\widehat\T} \widehat\eta^2( T)$).
Error estimator stability has been observed \cite{CFPP:AxiomsAdaptivity2014,CKNS:QuasioptimalConvergenceRate2008} for many
estimators and asserts in the 2-level notation at hand that
\begin{align}
	| \widehat\eta(\T\cap\widehat\T)-\eta(\T\cap\widehat\T)|\le\Lambda_1	\label{eqn:A1}\tag{A1}
\delta(\T,\widehat\T).
\end{align}
The slightly technical proof of stability \eqref{eqn:A1} utilizes reverse triangle inequalities and is well-established by now. 
For instance, the corresponding analysis in \cite{CH:HierarchicalArgyrisFinite2021}
is completely omitted. 
We therefore give an outline of the arguments and
emphasise the nonlinearity.
Indeed, the term $\eta(T)$ for $T\in\T$ is the Euclidean norm of a vector with the entries
$ |T|\|\lh\uh - \BiL \uh\|_{L^2(T)}$
as well as
$|T|^{1/4}\|\JF{\Dnn u_h}\|_{\LFace}$
and
$ |T|^{3/4}\|\JF{\DLn {u_h}} \|_{\LFace}$
for all $E\in\E(T)\cap\E(\GCSb)$.
The same interpretation of $\etat(T)$ for $T\in\Tt$ considers
$\eta(\T\cap \Tt) $ and $\etat(\T\cap\Tt)$ as Euclidean norms of vectors of the same length $m\in\N_0$.
Hence the squares of the reverse triangle inequality in $\R^m$ result in
\begin{align*}
 \left|  \widehat\eta(\T\cap\widehat\T)-\eta(\T\cap\widehat\T) \right|^ 2&\le 
 \sum_{T\in \T\cap\widehat\T} |T|^2\left| \| \lh\uh - \BiL \uh\|_{L^2(T)} -
\| \widehat\lh\widehat\uh - \BiL \widehat\uh\|_{L^2(T)} \right| ^2\\
&+ \sum_{T\in \T\cap\widehat\T} 
|T|^{1/2} \sum_{E\in\E(T)\setminus\E(\GCb)}
 \left| \|\JF{\Dnn {\widehat u_h}} \|_{\LFace}-\|\JF{\Dnn u_h}\|_{\LFace} \right| ^2 \\
&+ \sum_{T\in \T\cap\widehat\T} 
 |T|^{3/2} \sum_{E\in\E(T)\setminus\E(\GCSb)}   \left| 
\|\JF{\DLn{\widehat u_h}} \|_{\LFace}-\|\JF{\DLn { u_h}} \|_{\LFace} \right| ^2.
\end{align*}
 The reverse triangle
inequality in $L^2(T)$ resp.\ $L^ 2(E)$ allows for the upper bound
\begin{align}\nonumber
	\big|  \widehat\eta(\T\cap\widehat\T)&-\eta(\T\cap\widehat\T) \big|^ 2\le 
 \sum_{T\in \T\cap\widehat\T} |T|^2 \| (\widehat\lh\widehat\uh-\lh\uh)
 - \BiL\widehat e\|_{L^2(T)}^2\\
&+ \sum_{T\in \T\cap\widehat\T} 
|T|^{1/2} \sum_{E\in\E(T)\setminus\E(\GCb)}
 \|\JF{\Dnn {\widehat e}}\|_{\LFace}^2
+ \sum^{}_{T\in\T\cap\Tt} 
 |T|^{3/2} \sum_{E\in\E(T)\setminus\E(\GCSb)}   
\|\JF{\DLn{\widehat e}} \|_{\LFace}^2.\label{eqn:A1_reverse2}
\end{align}
The three sums in the upper bound are treated separately
and the non-linear volume terms
\[
 |T|^2 \| (\widehat\lh\widehat\uh-\lh\uh) - \BiL\widehat e\|_{L^2(T)}^2
 \le 2 |T|^2 \|\widehat\lh\widehat\uh-\lh\uh\|_{L^2(T)}^ 2
 +  2 |T|^2 \| \BiL\widehat e\|_{L^2(T)}^ 2
\]
in the first sum over $T\in \T\cap\widehat\T$
deserve special attention.
An inverse estimate in the last term leads to
some positive 
constant $C_\textrm{inv}$ (which exclusively depends on the shape-regularity in $\TT$
and the polynomial degree in the finite element space $V(\T)$) in 
\[
 |T| \| \BiL\widehat e\|_{L^2(T)}\le C_\textrm{inv} | \widehat e|_{H^2(T)}.
\]
The nonlinear term is estimated in \eqref{eqccadd13} and reveals 
\[
2\sum_{T\in \T\cap\widehat\T} |T|^2 \|\widehat\lh\widehat\uh-\lh\uh\|_{L^2(T)}^ 2
\le 	h_{\max}^4/2 \| \lht\uht  -\lh\uh\|^2  \le 
 h_{\max}^4  \lambda_{j+1} (2+\kappa)/8 \trb{ \widehat e}^2.
\]
The combination of the previous two displayed estimates leads to 
\[
\sum_{T\in \T\cap\widehat\T} |T|^2 
\| (\widehat\lh\widehat\uh-\lh\uh) - \BiL\widehat e\|_{L^2(T)}^2
 \le 
(h_{\max}^4  \lambda_{j+1} (2+\kappa)/8+2 C_\textrm{inv}^ 2)
 \trb{ \widehat e}^2.
\]
The routine estimation by triangle, trace, and inverse
inequalities as in \cite{CH:HierarchicalArgyrisFinite2021,Gra:OptimalMultilevelAdaptive2022} for the source problem
eventually bounds the jump terms in \eqref{eqn:A1_reverse2} by $\Lambda_1'^2\trb{\et}^2$.
This and the previous displayed estimate in \eqref{eqn:A1_reverse2} conclude the proof of \eqref{eqn:A1}
with
$\Lambda_1:=\sqrt{ {\Lambda_1' }^2+
\delta_2^4  \lambda_{j+1} (2+\kappa)/8+2 C_\textrm{inv}^ 2}$.\qed

\subsection{Proof of reduction (A2)}
Adopt the setting of the previous subsection and consider the set $\T\setminus\Tt$ of coarse but not fine
triangles together with the fine but not coarse triangles in $\widehat\T\setminus\T$.
The sum convention defines the
two related estimator terms $\eta(\T\setminus\widehat\T)$ resp.\
$\widehat\eta(\widehat\T\setminus\T)$ on the coarse resp.\ fine level. 
The aim of this subsection is the proof of the error estimator
reduction 
\begin{align}
	\widehat\eta(\widehat\T\setminus\T)\le 2^{-1/4} 
	\eta(\T\setminus\widehat\T)+\Lambda_2\delta(\T,\widehat\T).\label{eqn:A2}\tag{A2}
\end{align}
The main arguments are those from the previous section and in fact 
$\Lambda_2=\Lambda_1$. But the details are more involved and perhaps simplified with
auxiliary estimator terms $\widetilde\eta^2(T) $ defined for each 
fine triangle $T\in\widehat\T$ with respect to the coarse eigenpair 
$(\lh,\uh)\in\R^ +\times V(\T) $, i.e.,
\begin{align}\label{eqn:eta_deffinetilde}
\widetilde\eta^2( T) &= |T|^2\| \lh\uh - \BiL\uh\|_{L^2(T)}^2\\ 
\nonumber
&+|T|^{1/2} \sum_{E\in\Et(T)\setminus\widehat\E(\GCb)} 
\|\JF{\Dnn  u_h}\|_{\LFace}^2
+  |T|^{3/2}\sum_{E\in\Et(T)\setminus\widehat\E(\GCSb)}  \nonumber
\|\JF{\DLn { u_h}} \|_{\LFace}^2.
\end{align}
Notice the same index sets over the triangle and edges with respect to the fine
geometry of $\Tt$
in \eqref{eqn:eta_deffine} and \eqref{eqn:eta_deffinetilde}. 
Reverse triangle
inequalities lead 
(in $\R^m$ and in $L^2(T)$ resp.\ $L^ 2(E)$)
as in the previous subsection~to
\begin{align*}
\widehat\eta(\widehat\T\setminus\T)
&-\widetilde\eta(\widehat\T\setminus\T)
 \le 
 \sum_{T\in \widehat\T\setminus\T} |T|^2 \| (\widehat\lh\widehat\uh-\lh\uh)
 - \BiL\widehat e\|_{L^2(T)}^2\\
&+ \sum_{T\in \widehat\T\setminus\T} 
|T|^{1/2} \sum_{E\in\Et(T)\setminus\widehat\E(\GCb)}
 \|\JF{\Dnn {\widehat e}}\|_{\LFace}^2 
+ \sum_{T\in \widehat\T\setminus\T} 
 |T|^{3/2} \sum_{E\in\Et(T)\setminus\widehat\E(\GCSb)}   
 \|\JF{\DLn{\widehat e}} \|_{\LFace}^2.
\end{align*}
The summands in the upper bound are identical to the situation in the previous subsection and the same argumentation
applied on the fine level for $\Tt\setminus\T$ leads to
\begin{align}\label{eqn:A2_intermediate}
\widehat\eta(\widehat\T\setminus\T)
\le \widetilde\eta(\widehat\T\setminus\T)
+\Lambda_2\delta(\T,\widetilde\T).
\end{align}
The proof of the remaining
reduction property $ \widetilde\eta(\widehat\T\setminus\T)\le 2^{-1/4} 
\eta(\T\setminus\widehat\T)$ rewrites the sum over all fine but not coarse triangles 
$T\in \widehat\T\setminus\T$ as a double sum over all coarse but not fine triangles
$K\in \T\setminus\widehat\T$ and its children in 
$\widehat\T(K):=\{T\in\widehat\T: T\subset K\}$.
The point is that, for each $K\in \T\setminus\widehat\T$, a triangle $T\in\widehat\T(K)$ has been obtained by bisection
of $K$ at least once so that 
$|T|\le |K|/2$ implies
\begin{align*}%
	&\sum_{T\in \widehat\T(K)}\widetilde \eta^2( T) \le  
2^ {-2} |K|^2\sum_{T\in \widehat\T(K)} \| \lh\uh - \BiL\uh\|_{L^2(T)}^2\\ 
\nonumber
&\quad+2^ {-1/2}|K|^{1/2} \sum_{T\in \widehat\T(K)} 
\bigg(\sum_{E\in\E(T)\setminus\widehat\E(\GCb)} 
\|\JF{\Dnn  u_h}\|_{\LFace}^2
+  2^{-1}|K|
\sum_{E\in\E(T)\setminus\widehat\E(\GCSb)}  \nonumber
\|\JF{\DLn { u_h}} \|_{\LFace}^2\bigg).
\end{align*} 
The sums over $T\in \widehat\T(K)$ combine directly in 
$\sum_{T\in \widehat\T(K)} \| \lh\uh - \BiL\uh\|_{L^2(T)}^2=\| \lh\uh - \BiL\uh\|_{L^2(K)}^2$,
while the jumps with respect to the fine geometry vanish for interior edges 
$\widehat E\in \widehat\E$ with $\widehat E\not\subset\partial K$ for  $u_h$
is a polynomial in $K$. Indeed, the non-zero summands for 
$\widehat E\in\E(T)$, $T\in \widehat\T(K)$,
belong to the boundary  $\widehat E\subset\partial K$ and sum up 
to the full $L^2$ norms along the coarse edges $E$ of $K$. This explains  
\begin{align*}
 \sum_{T\in \widehat\T(K)} 
\sum_{E\in\E(T)\setminus\widehat\E(\GCb)} 
\|\JF{\Dnn  u_h}\|_{\LFace}^2&= \|\JF{\Dnn  u_h}\|_{L^ 2(\partial K\setminus \GCb)}^2,
\\
\sum_{T\in \widehat\T(K)}
\sum_{E\in\E(T)\setminus\widehat\E(\GCSb)}  \nonumber
\|\JF{\DLn { u_h}} \|_{\LFace}^2&=\|\JF{\Dn\Delta  u_h}\|_{L^ 2(\partial K\setminus \GCSb)}^2.
\end{align*} 
The combination of the above estimates reveals (with the sum convention) that
\[
\widetilde \eta^2( \widehat\T\setminus\T)=\sum_{T\in \widehat\T(K)}\widetilde \eta^2( T)\le 2^{-1/2} \eta^
2(K)\qquad\text{for }K\in \T\setminus\widehat\T.
\]
This and \eqref{eqn:A2_intermediate} conclude the proof of \eqref{eqn:A2}.\qed

\subsection{Proof of discrete reliability (A3)}
This subsection provides the discrete reliability in the 2-level notation of 
Subsection~\ref{sub:2levelnotationConvergenceconformingEVP}: Given 
$\T\in\TT(\delta_3)$ with  $\delta_3\le\delta_2$ as in Lemma~\ref{prop:stab} and  
\ref{lem:err_Res_est}. In fact, those preliminaries form the key arguments for
\begin{align}
	\delta(\T,\widehat\T)=\trb{ \widehat e}\le \Lambda_3\eta(\T\setminus\widehat\T).\label{eqn:A3}\tag{A3}
\end{align}
The proof of (A3)  is finalized by the known arguments on \cite{CH:HierarchicalArgyrisFinite2021} for the associated source problem
with $f:=\lh \uh\in L^ 2(\Omega)$. 
In fact,  Lemma~\ref{lem:err_Res_est} provides 
\[
(1-\kappa) \trb{\et}^2\le\textrm{Res}(\et)=\int_\Omega f\,\et\, \mathrm dx-a(\uh,\et)
\]
and the latter residual is already considered in \cite[Thm.~4]{CH:HierarchicalArgyrisFinite2021}. With an identical proof,
based on a discrete quasi-interpolation  operator \cite[Thm.~2]{CH:HierarchicalArgyrisFinite2021} for $p=5$
and further nowadays standard arguments in explicit 
residual-based a posteriori error control,  that former result reveals 
\[
\int_\Omega f\,\et\, \mathrm dx-a(\uh,\et)\le \Lambda_3' \trb{\et} \eta(\T\setminus\widehat\T).
\]
The universal constant $\Lambda_3'$ (called $\Lambda_3$ in
\cite{CH:HierarchicalArgyrisFinite2021,Gra:OptimalMultilevelAdaptive2022})
depends on the polynomial degree of the Argyris
finite element scheme and the shape regularity in $\TT$; 
 $\Lambda_3'$  is independent of the discrete
or continuous eigenvalues (i.e. independent of $j$ fixed in the 2-level notation) and
independent of any smallness assumption (i.e. independent of $\delta_3$).
This concludes the proof of (A3) with $\Lambda_3:= \Lambda_3'/(1-\kappa)$ for $p=5$. 
The generalisation of \cite[Thm.~2]{CH:HierarchicalArgyrisFinite2021} for $p\ge 6$ is possible with the techniques of
\cite{CH:HierarchicalArgyrisFinite2021}.\qed

\subsection{Proof of quasi-orthogonality  (A4)}
This is in fact an orthogonality based on the nestedness of the
discrete spaces $V(\T_\ell)\subset V(\T_{\ell+1})$ already utilized in 
 Lemma~\ref{prop:stab} and   \ref{lem:err_Res_est}. The estimate
 \(
(1-\kappa) \trb{\et}^2\le \lh-\lht
\)
of  Lemma~\ref{lem:err_Res_est}.a  translates for 
$\T:=\T_\ell $ and 
$\widehat\T:=\T_{\ell+1}$ into
\[
(1-\kappa)\delta^ 2(\T_\ell, \T_{\ell+1}) \le \lambda_\ell - \lambda_{\ell+1}
\]
in terms of the output of the adaptive algorithm. For any $\ell,m\in\N_0$,
this implies a telescoping sum
\begin{align*}
(1-\kappa)	\sum^{\ell+m}_{k = \ell} 
\delta^ 2(\T_\ell, \T_{\ell+1})
\leq \lambda_\ell - \lambda_{\ell+m+1}.
\end{align*}
Lemma~\ref{lem:err_Res_est}.a and (A3) for to $\T:=\T_\ell $ and 
$\widehat\T:=\T_{\ell+m+1}$ provide
\[
 \lambda_\ell - \lambda_{\ell+m+1}\le \delta^ 2(\T_\ell, \T_{\ell+m+1})\le \Lambda_3^2
 \eta^2(\T_\ell).
\]
The combination of the last two displayed estimates leads to
$\Lambda_4:= \Lambda_3/(1-\kappa)$ in 
\begin{align}
\sum^{\ell+m}_{k = \ell} 
\delta^ 2(\T_\ell, \T_{\ell+1})\le \Lambda_4  \eta^2(\T_\ell)
\quad\text{for all }\ell,m\in\N_0.
\label{eqn:A4}\tag{A4}\qed
\end{align}

\subsection{Proof of optimal rates}\label{sub:proof_optimal_rates}
This paper has established  \eqref{eqn:A1}--\eqref{eqn:A3}
 for any triangulation $\T\in\TT(\delta_2)$ and its refinement  $\T\in\TT(\T)$ and 
 \eqref{eqn:A4} for the outcome of  AFEM provided $\T_0\in\TT(\delta_2)$.
 It is well established, that those properties imply optimality in the sense of \eqref{eqccmainresultinpaper} and
 conclude the proof of Theorem~\ref{thmoptimalrates}.
 This is proven explicitly in 
 \cite{CFPP:AxiomsAdaptivity2014, CR:AxiomsAdaptivitySeparate2017}
 and also follows with arguments from 
  \cite{CKNS:QuasioptimalConvergenceRate2008,Ste:CompletionLocallyRefined2008}.\qed
  
\FloatBarrier
\section{Numerical results}%
\label{sec:Numerical Results}
This section provides striking numerical evidence of the optimal
convergence rates for the proposed hierarchical adaptive Argyris FEM and thereby displays high-precision reference eigenvalues.

\subsection{Comments on the implementation}%
\label{sub:Numerical realisation}
Analytic exact eigenvalues are \emph{not} available
for basic geometries including the
square domain in Subsection~\ref{sub:high_precision}.
The rapid convergence of the high-order Argyris AFEM %
results in a new state-of-the art method for the computation of
reference eigenvalues at the expense of a large condition number
of the algebraic
eigenvalue problem.
Our realisation of the Argyris AFEM described in \cite{Gra:OptimalMultilevelAdaptive2022} therefore
utilizes multi-precision arithmetic 
\cite{Zha:MultiFloatsJl2023}. %
We evaluate integrals with an exact quadrature rule and expect all displayed digits in this section to be significant.
An implicitly
restarted Arnoldi method for the algebraic eigenvalue problem $Ax=\lambda Bx$ with mass matrix $B$ and stiffness
matrix $A$ using the shift-and-invert strategy~\cite{Sco:AdvantagesInvertedOperators1982,GLS:ShiftedBlockLanczos1994} solves the shifted problem $(A-\sigma B)^{-1}B x = \alpha x$ with
eigenvalues $\lambda=\sigma+\alpha^{-1}$ for the shift $\sigma\geq0$ with $\sigma=0$ if not stated otherwise.
\W{We note that the efficient solution of the algebraic eigenvalue problem, e.g., by preconditioned multilevel schemes available for
the source
problem~\cite{BZ:MultigridMethodsBiharmonic1995,AP:PreconditioningHighOrder2021,CH:HierarchicalArgyrisFinite2021,Gra:OptimalMultilevelAdaptive2022},
is left for future research.}

This section considers clamped boundary conditions ($u|_{\Gamma_C}\equiv0\equiv u|_{\Gamma_C}$) on a closed part $\Gamma_C\subset\partial\Omega$ and simply-supported boundary
conditions ($u|_{\Gamma_S}\equiv 0$) on a relatively open part $\Gamma_S\subset\partial\Omega$ of the boundary.
The hierarchical Argyris AFEM for the mixed boundary conditions is discussed in \cite{Gra:OptimalMultilevelAdaptive2022}
and leads to the modified error estimator $\eta=\big(\sum^{}_{T\in\T} \eta^2(T)\big)^2$ with the local contributions
\begin{align*}
	\eta^2( T) &= |T|^2\|\lh\uh - \BiL \uh\|_{L^2(T)}^2\\
			   &+|T|^{1/2}\|\JF{\partial_{\nu\nu}^2  u_h}\|_{L^2(\partial T\setminus \Gamma_C)}^2+  |T|^{3/2}
\|\JF{\partial_{\tau\tau\nu}^3u_h+\partial_{\nu}\Delta u_h}\|_{L^2(\partial T\setminus(\Gamma_C\cup\Gamma_S))}^2.
\end{align*}
The differences of this estimator to \eqref{eqn:eta_def} are the modified jumps on the free boundary
$\Gamma_F\coloneqq\partial\Omega\setminus(\Gamma_C\cup \Gamma_S)$.
On interior edges $E\in\E(\Omega)$, $\JF{\partial_{\tau\tau\nu}^3u_h}$ vanishes by the $C^1$ continuity of the discrete
solution $u_h\in V(\T)$ and the jumps coincide with those in \eqref{eqn:eta_def}.
The optimal convergence rates from Section~\ref{sec:Argyris FEM for eigenvalue problems} generalizes to the
mixed boundary conditions of this problem. \W{The benchmarks compare AFEM (with $\theta=0.5$ unless stated otherwise)
	and uniform red-refinement, where each
triangle is divided into four congruent triangles by joining the edge midpoints}; ndof abbreviates the number of degrees
of freedom.

\begin{figure}[]
	\centering
	\hspace*{-4em}\hbox{
	\includegraphics[width=.28\textwidth]{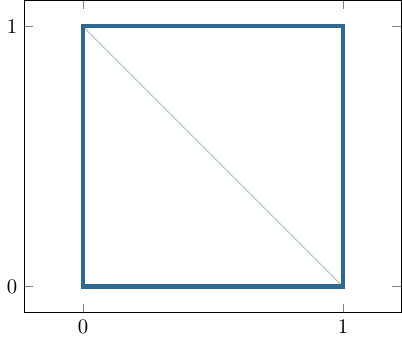}
	\includegraphics[width=.28\textwidth]{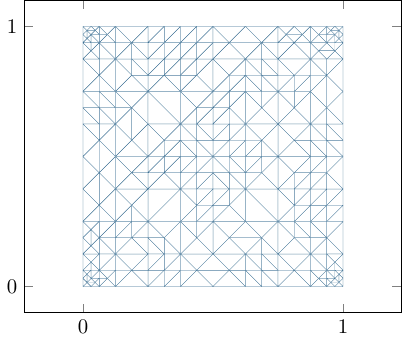}
	\includegraphics[width=.28\textwidth]{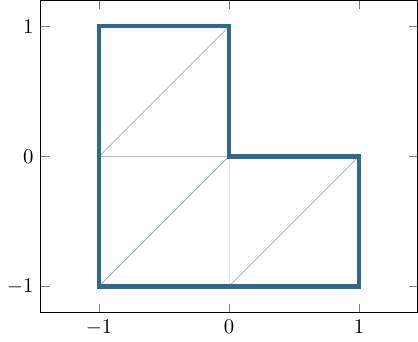}
	\includegraphics[width=.28\textwidth]{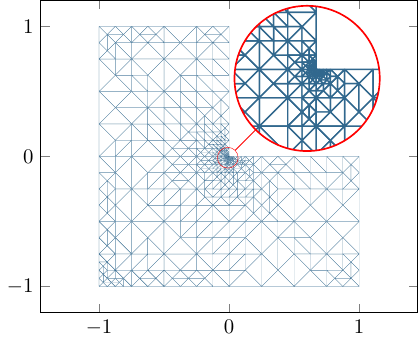}
}
	\caption{Initial and adapative triangulations of the unit square and the L-shaped domain}
	\label{fig:Square_Lshape_mesh}
\end{figure}
\subsection{High-precision eigenvalues on the square}%
\begin{figure}[ht]
		\centering
			\includegraphics{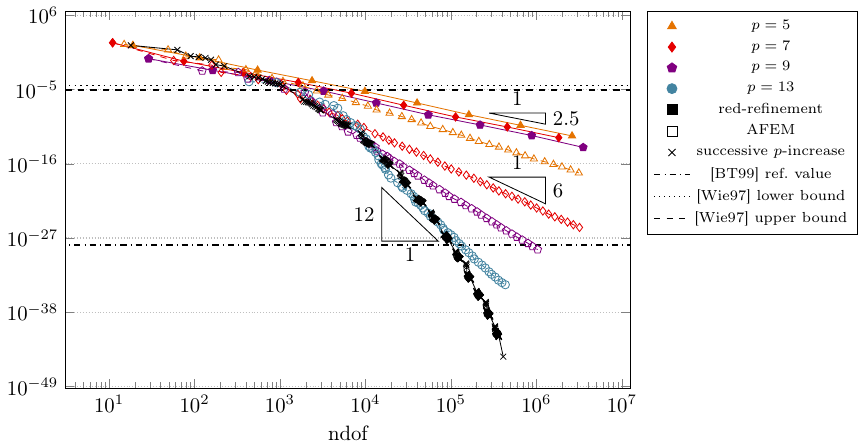}
			\caption{Convergence history plot of AFEM \W{(hollow) and uniform-red refinement (filled)} towards the principal eigenvalue $\lambda_1$ of the unit square for different $p$ and a successive
		$p$-increase}%

				\label{fig:Square_h_hp}
	\end{figure}

\label{sub:high_precision}
The approximation of the principal eigenvalue on the unit square $\Omega=(0,1)^2$ results in the reduced
empirical convergence rate $2.5$ in terms of ndof under uniform mesh refinement shown in
Figure~\ref{fig:Square_h_hp}.
\begin{table}[]
	\centering
	\begin{tabular}{rl}
		&Square\\
		\hline
	$\lambda_1$&				$1294.9339795917128081703026479743085522513148$ \\%5$\\
	$\lambda_2=\lambda_3$&		$5386.6565607779451709440883164319500534323747$ \\%4$ \\
	$\lambda_4$&				$11710.811238205718716479524026744165110548790$ \\
	$\lambda_5$&				$17313.499721776700784267277477409032730611369$ \\%455$\\
	$\lambda_6$&				$17478.106551478646880691169158961070602062916$\\
	$\lambda_7=\lambda_8$&		$27225.134042723468407942851261240075600133838$\\
	$\lambda_9=\lambda_{10}$&	$44319.444997239644605551889845792194008818266$\\
	\end{tabular}
	\caption{First ten eigenvalues of the unit square}
	\label{tab:ref_val_square}
\end{table}

AFEM for $\theta=0.5$ and $p=13$ computes the reference values in Table~\ref{tab:ref_val_square} on an
adaptive mesh with $2\times10^6$ ndofs in between $70$--$85$ levels.
The multiple eigenvalues $\lambda_2=\lambda_3, \lambda_7=\lambda_8$, and $\lambda_9=\lambda_{10}$ lead to numerical
instability of the Arnoldi method if no shift $\sigma=0$ was applied.
In this case, only up to 32 significant digits of $\lambda_3$ to $\lambda_{10}$, independent of the precision in the
computations, could be obtained.
The remedy with a shift leads to the reference values in Table~\ref{tab:ref_val_square}: 
the shift $\sigma$ therein is
the integral part of the respective eigenvalue.
Figure~\ref{fig:Square_h_hp} displays the optimal convergence rates of the Argyris AFEM for different polynomial
degrees $p$ that are guaranteed by Theorem~\ref{thmoptimalrates} also for
a modification of AFEM with $\theta=0.1$ and a successive $p$-increase.
In the successive $p$-increase, the AFEM remains as described in Figure~\ref{fig:AFEM} but for variable polynomials
degrees: For $\ell=0$, let $p=5$ and increase $p\coloneqq p+1$ after additional $2p$ levels, i.e.,
$p=5$ for $\ell=0,\dots,9$, $p=6$ for $\ell=10,\dots,21$, $p=7$ for $\ell=22,\dots,35$, etc.
The latter strategy offers a highly efficient alternative with
super-algebraic convergence rates at comparable computational costs.
The horizontal lines in Figure~\ref{fig:Square_h_hp} indicate the known eigenvalue inclusions from
\cite{Wie:BoundsLowestEigenvalues1997} and the high-precision reference value from \cite{BT:HighPrecisionSolutions1999}.

	\FloatBarrier
\subsection{L-shaped domain}%
\begin{figure}[ht]
		\centering
			\hbox{\hspace{-5em}\scalebox{0.9}{\includegraphics{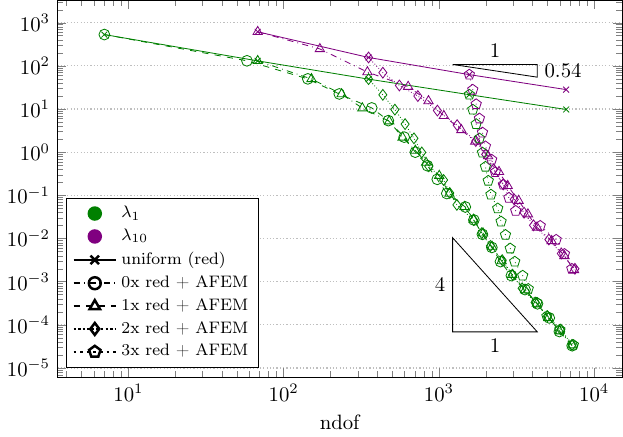}}
			\scalebox{0.9}{\includegraphics{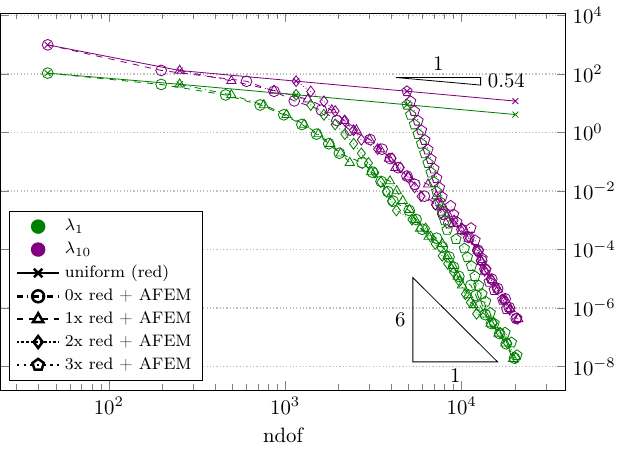}}
			}
			\caption{Convergence history plot of $\lambda_1$ and $\lambda_{10}$ from AFEM ($p=5$ left and $p=7$
				right) with $n=0,1,2,3$-times red-refinement of $\T_0$ from
			Figure~\ref{fig:Square_Lshape_mesh} as initial triangulation}

		\label{fig:LshapeConvergenceAB}
	\end{figure}

\label{sub:Experiments on the L-shaped domain}
The L-shaped domain $\Omega\coloneqq (-1,1)^2\setminus[0,1)^2$ with initial triangulation in
Figure~\ref{fig:Square_Lshape_mesh} has an index of elliptic regularity $\sigma_{\textrm{reg}}=0.54$
\cite{Gri:SingularitiesBoundaryValue1992}.
Figure~\ref{fig:LshapeConvergenceAB} displays 
the expected suboptimal convergence rate $\sigma_{\textrm{reg}}=0.54$ (in terms of the degrees of freedom) under uniform refinement for $\lambda_1$ and
$\lambda_{10}$.
The adaptive algorithm refines towards the singularity at the origin $(0,0)$ as shown in Figure~\ref{fig:Square_Lshape_mesh} for
$p=5$ and recovers optimal convergence rates $p-1$ asymptotically as predicted by Theorem~\ref{thmoptimalrates} for a
sufficiently small maximal mesh-size of the initial triangulation.

AFEM recovers optimal asymptotic rates already for very coarse initial triangulations in
Figure~\ref{fig:LshapeConvergenceAB}, but
exhibits a pre-asymptotic behaviour with reduced convergence rates on coarse meshes.
Initial red-refinements offer a cheap way to reduce the computational time in the pre-asymptotic regime
of suboptimal convergence which length depends on the polynomial degree
$p$ and the approximated eigenvalue.
Figure~\ref{fig:LshapeConvergenceAB} compares AFEM for $\theta=0.5$ and $p=5$ applied to $n=0,1,2,3$ (displayed with differently shaped markers) initial red-refinements of the initial
triangulation from Figure~\ref{fig:Square_Lshape_mesh}.  
Already two initial red-refinements remove the pre-asymptotic regime almost completely, whereas 
an even higher number of initial red-refinements lead to a superconvergence behaviour of AFEM \W{with many small
increments of ndof} until the convergence graphs meet and exhibit the same convergence error. 

Table~\ref{tab:ref_val_L} provides reference values of the first ten eigenvalues on the L-shaped domain.

\begin{table}[]
	\centering
	\begin{tabular}{rl}
		&L-shape\\
		\hline
		$\lambda_1$&				$418.9752928519954616775618775449403$\\
		$\lambda_2$&				$690.9065117037020521173375432095840$ \\
		$\lambda_3$&				$931.5792655819233500496763905727078$ \\
		$\lambda_4$&				$1634.533781725410909450276640848249$\\
		$\lambda_5$&				$2090.839376830117972780944169183919$\\
		$\lambda_6$&				$3350.410627882138361796676879507818$\\
		$\lambda_7$&				$3720.925878799246215492176662379798$\\
		$\lambda_8$&				$4485.620042035762239796922517307530$\\
		$\lambda_9$&				$4560.432156327978825966680488256231$\\
		$\lambda_{10}$&				$5738.403842705743348891002113406246$\\
	\end{tabular}
	\caption{Reference values for the L-shaped domain}
	\label{tab:ref_val_L}
\end{table}

	\FloatBarrier

\subsection{Isospectral domains}%
\begin{figure}[ht]
	\centering
	\scalebox{0.6}{\includegraphics{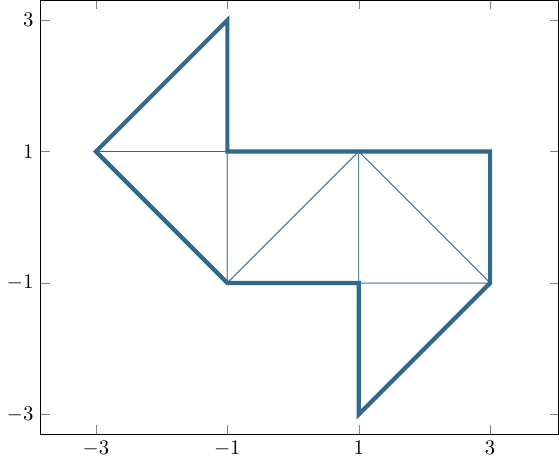}}
	\scalebox{0.6}{\includegraphics{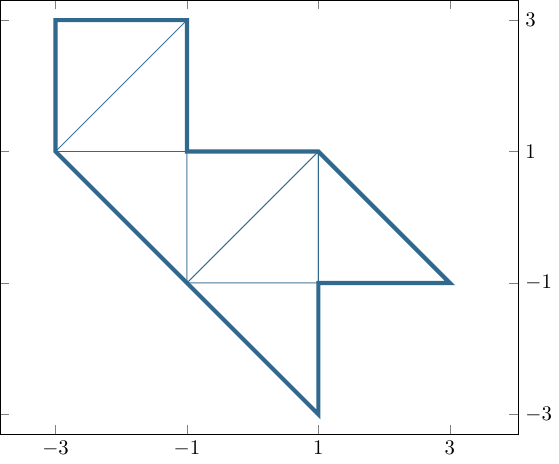}}
	\caption{Two isospectral domains for the Laplacian or the biharmonic operator under simply-supported boundary
	conditions}
	\label{fig:Drum_mesh}
\end{figure}

\label{sub:Experiments on isospectral}
The spectrum of the Laplacian for different
domains can coincide
\cite{GWW:OneCannotHear1992,Bus:IsospectralRiemannSurfaces1986,BCDS:PlanarIsospectralDomains1994,Kac:CanOneHear1966}; 
Figure~\ref{fig:Drum_mesh} displays the initial triangulation of two Laplace-isospecral domains.
Our computations provide numerical evidence that these domains are also isospectral for the biharmonic operator under
simply-supported boundary conditions
($u=0$ on
$\partial\Omega$ in $V=H^2(\Omega)\cap H^1_0(\Omega)$), cf.~Remark~\ref{rem:isospectral} for an explanation.
Table~\ref{tab:ref_Drum} provides reference values of the biharmonic operator for the isospectral domains of
Figure~\ref{fig:Drum_mesh} under simply-supported boundary conditions and verifies that these domains are \emph{not}
isospectral under clamped boundary conditions ($u=0=\partial_\nu u$ on $\partial\Omega$) of \eqref{eqn:WP}.
The ninth eigenfunction $u_9(x,y)=\cos(\pi x/2)\sin(\pi y) - \cos(\pi y/2)\sin(\pi x)\in C^\infty(\Omega)\cap H^1_0(\Omega)$ of
the Dirichlet-Laplacian on both domains in Figure~\ref{fig:Drum_mesh} is
special~\cite{Dri:EigenmodesIsospectralDrums1997}.
Restricted to each right-isosceles triangle $T\in\T$ in Figure~\ref{fig:Drum_mesh}, $u_9|_T$ is the first eigenmode of the
Laplacian (on that
triangle) described in~\cite{DPCD:UnderstandingEigenstructureVarious2010} with eigenvalue $5\pi^2/4$.
Since 
$u_9\in H^2(\Omega)\cap H^1_0(\Omega)$,
it is also an eigenfunction of the biharmonic operator with simply-supported boundary
conditions for the eigenvalue $(5\pi^2/4)^2=\lambda_9$ in Table~\ref{tab:ref_Drum}.
\begin{figure}[ht]
	\centering
		\hbox{\hspace{-5em}\scalebox{0.9}{\includegraphics{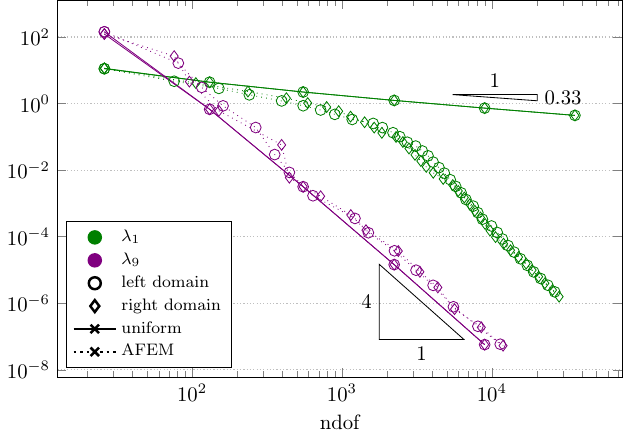}}
		\scalebox{0.9}{\includegraphics{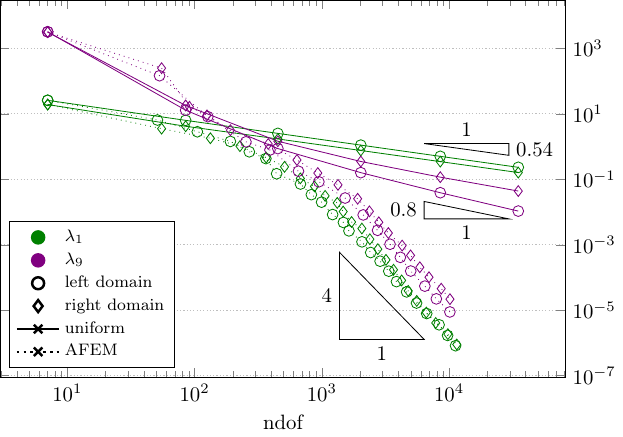}}
		}
		\caption{Convergence history plot of $\lambda_1$ and $\lambda_{9}$ on the two domains of
			Figure~\ref{fig:Drum_mesh} with simply-supported (left) and clamped boundary conditions (right)}

	\label{fig:IsospectralConvergenceAB}
\end{figure}

The smoothness of $u_9$ leads in
Figure~\ref{fig:IsospectralConvergenceAB} to optimal convergence rates under both uniform and adaptive
refinement for $p=5$.
This is in contrast to the principal eigenvalue of the biharmonic under simply-supported
boundary conditions with a small convergence of rate $0.33$ under uniform refinement.
No localisation of eigenfunctions over the congruent
triangles of Figure~\ref{fig:Drum_mesh} similar to that of $u_9$ has been observed for the biharmonic operator with clamped
boundary conditions.
AFEM recovers optimal convergence rates
in all cases.
The observed convergence rate $0.8$ (on the right) for the ninth eigenfunction under uniform refinement over the displayed
computational domain
may be preasymptotic.

\begin{remark}[Isospectrality for simply-supported boundary conditions]\label{rem:isospectral}
The proof of isospectrality for two domains $\Omega_1$ and $\Omega_2$ by the transplantation technique from Buser
\cite{Bus:IsospectralRiemannSurfaces1986}
for the Laplace operator
explicitly
constructs a map $\Phi:H^1_0(\Omega_1)\to H^1_0(\Omega_2)$
that transforms an eigenfunction $u_1(\lambda)$ of $\Omega_1$ to
an eigenfunction $u_2(\lambda)$ with the same eigenvalue $\lambda$ on the second domain $\Omega_2$.
In short, that technique relies on triangulations $\T_1$ and $\T_2$ of $\Omega_{1}$ and $\Omega_2$ into congruent
triangles as in
Figure~\ref{fig:Drum_mesh}.
The idea is then to locally construct an eigenfunction 
$u_2(\lambda)|_K$ on a triangle $K\in\T_2$ in $\Omega_2$ for the same eigenvalue $\lambda$ as a linear combination of
re-transformed $\{u_1(\lambda)|_T\}_{T\in
\T_1}$.
Then $u_2\in H^1_0(\Omega)$ becomes an eigenfunction on $\Omega_2$
with eigenvalue $\lambda$.
A list of domain pairs that allow such a design is studied in
\cite{BCDS:PlanarIsospectralDomains1994}.
A careful investigation of the continuity and boundary
conditions of the transplanted function verifies that this technique also applies to the biharmonic operator
under simply-supported
boundary conditions but \emph{fails} for the clamped boundary conditions.
This justifies the experimentally verified isospectrality only for simply-supported boundary conditions, but fails
for clamped boundary conditions.

\end{remark}

\newcommand{\lDa}{\lambda}
\newcommand{\lDb}{\lambda^{(2)}}
\begin{table}[]
	\centering
	\hspace*{-2em}\begin{tabular}{rllll}
		&\hspace*{-1em}simply-supported (both domains)& clamped (left domain)&clamped (right domain)\\
		\hline
	$\lDa_1$&$10.36402107986949540721$	&$28.0586863865423000549087398$		&	$25.01410502064436259175268775$		\\
	$\lDa_2$&$16.402184923407268356740$	&$42.80755796191819473005366012$	&	$50.67098031228146895086815044$		\\
	$\lDa_3$&$36.28192141809921686643$	&$73.51576421441210956574548706$	&	$72.0887124830690440717331005$		\\
	$\lDa_4$&$47.33781364226172710396$	&$109.355356090831716811598988$		&	$101.50786786148196218598065170$		\\
	$\lDa_5$&$61.23040947172451684072$	&$123.8961882657907045115913099$	&	$119.568610048818932875544158$		\\
	$\lDa_6$&$90.01373420512098390986$	&$151.830513099492311414264541$		&	$162.5851723649249044594031075$	\\
	$\lDa_7$&$119.052506025401956876870$	&$227.038626156347910723847312$		&	$217.9878528191822288137668706$	\\
	$\lDa_8$&$135.11181644660087949787$	&$251.63057739930330617446625802$	&	$263.6684022735138072395955186$	\\
	$\lDa_9$&$152.2017047406288081819380$&$290.252435400071762741097143187$	&	$297.5045897085804999794427456$	\\
 $\lDa_{10}$&$187.12311642208418946283$	&$311.16810084585353295577573668$	&	$315.5017873309864911264409052$	\\
	\end{tabular}
	\caption{Eigenvalues of the biharmonic operator for the two domains in Figure~\ref{fig:Drum_mesh} under
		simply-supported
		(left) and
		under clamped boundary conditions (middle and right)
	}
	\label{tab:ref_Drum}
\end{table}

	\FloatBarrier
\subsection{Experiments on the rectangle with hole}%
\begin{figure}[ht]
	\centering
	\scalebox{0.6}{\includegraphics{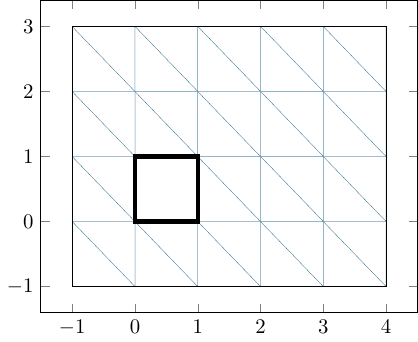}}
	\scalebox{0.6}{\includegraphics{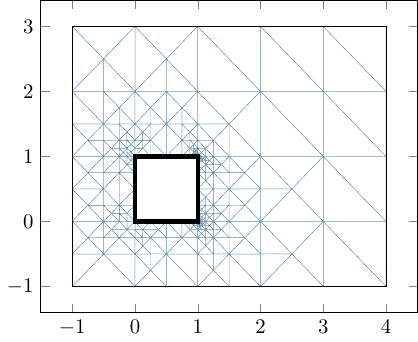}}
	\scalebox{0.6}{\includegraphics{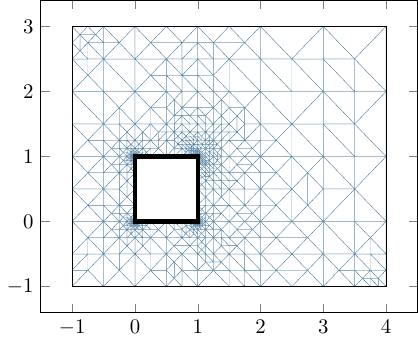}}
	\caption{Initial triangulation (left)
	and adaptive refinements with $|\T|=572$ (middle) and $|\T|=2528$ (right) for the principal eigenvalue and $p=5$ in
Subsection~\ref{sub:Experiments on a rectangle with a hole}}
	\label{fig:RectangleHole}
\end{figure}

\label{sub:Experiments on a rectangle with a hole}
This benchmark considers the domain
$\Omega\coloneqq R\setminus [0,1]^2$ displayed in Figure~\ref{fig:RectangleHole} 
with free boundary conditions on the outer boundary $\Gamma_F=\partial R$ of the
rectangle $R= (-1,3)\times(-1,4)$ and clamped boundary conditions on the interior boundary
$\Gamma_C=\partial[0,1]^2$. %
The abstract setting of Section~\ref{sec:A posteriori} applies for this problem with mixed boundary conditions 
in the space of admissible functions $V\coloneqq\{v\in H^2_0(\Omega)\ :\ v|_{\Gamma_C}=0=\partial_\nu
v|_{\Gamma_C}\}$.

\begin{figure}[ht]
	\centering
	\hspace*{-4em}\hbox{
	\includegraphics[]{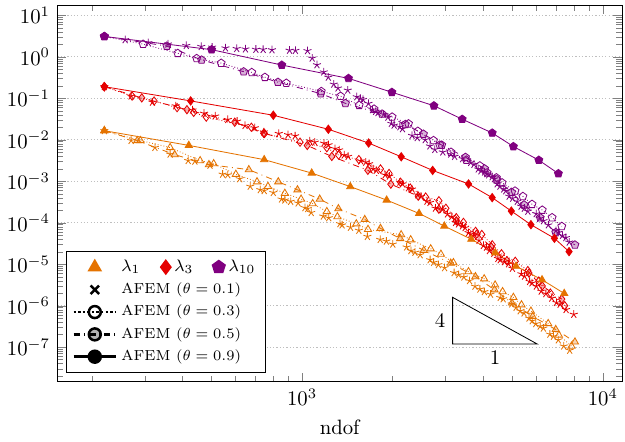}
	}
	\caption{Convergence history plot of AFEM with different bulk parameters $\theta$ for some eigenvalues
in Subsection~\ref{sub:Experiments on a rectangle with a hole}}
	\label{fig:RectangleHoleConvergence}
\end{figure}

Undisplayed computations reveal the empirical convergence rate $0.54$ (in terms of ndof) for the approximation of the principal
eigenvalue under uniform mesh-refinement.
This suggests that the corresponding eigenfunction has the same regularity as for the L-shaped domain of
Subsection~\ref{sub:Experiments on the L-shaped domain}.
AFEM for the initial triangulation of Figure~\ref{fig:RectangleHole} locally refines towards the four reentering corners
with different intensity (depending on the approximated eigenvalue) and recovers optimal convergence rates in
Figure~\ref{fig:RectangleHoleConvergence} for $p=5$ already on coarse triangulations.
A smaller bulk parameter $\theta$ results in a lower error on fine meshes for the approximation of three displayed
eigenvalues, at the cost of an increased number of iterations in the AFEM loop.
However, the asymptotic improvement of $\theta<0.5$ over $\theta=0.5$ appears negligible.
If the discrete spectrum
is not resolved sufficiently accurate, a small value 
$\theta=0.1$ for the approximation of $\lambda_{10}$ in
Figure~\ref{fig:RectangleHoleConvergence} even enlarges the preasymptotic regime.
One uniform red-refinement (undisplayed) of the initial mesh removes the observed preasymptotic behaviour and leads to the guaranteed optimal
rates from the first iteration with ndof=$786$.

\begin{table}[]
	\centering
	\begin{tabular}{rl}
		$\lambda_1$&	$0.10698498562334817102814013$\\
		$\lambda_2$&	$0.35605676603875088420438615$ \\
		$\lambda_3$&	$0.94524070807442103593431671$ \\
		$\lambda_4$&	$2.53208704546115218637535529$\\
		$\lambda_5$&	$3.99930885285784387075239865$\\
		$\lambda_6$&	$5.57646932761755697998226709$\\
		$\lambda_7$&	$6.43857726128711796898224149$\\
		$\lambda_8$&	$9.13840065843740557418337614$\\
		$\lambda_9$&	$16.4479467402360010159306274$\\
		$\lambda_{10}$&	$17.4245866515989760773466203$\\
	\end{tabular}
	\caption{First ten eigenvalues for the rectangle with a hole}
	\label{tab:ref_val_rectangle_hole}
\end{table}
\newpage
\subsection{Constants related to eigenvalues on triangles}%
\label{sub:Special constants on triangles}
\begin{figure}[ht]
	\centering
	\scalebox{0.6}{\includegraphics{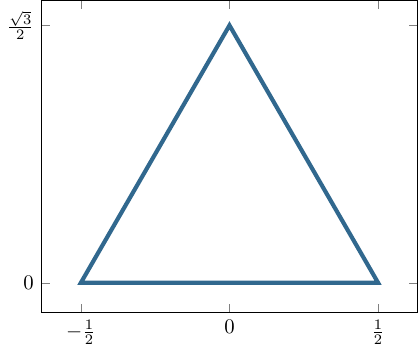}}
	\scalebox{0.6}{\includegraphics{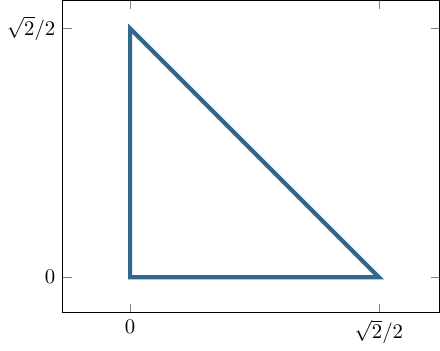}}
	\scalebox{0.6}{\includegraphics{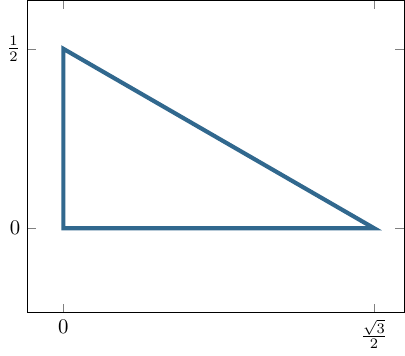}}
	\label{fig:Triangles}
	\caption{Equilateral, Right-isosceles, and 90-60-30 triangle scaled to $h_T=1$}
\end{figure}
This benchmark reports on some eigenvalues related to constants in standard interpolation estimates for the equilateral triangle, right-isosceles triangle, and the right triangle with angles
90-60-30 depicted in Figure~\ref{fig:Triangles}.
Let $T\subset \R^2$ be a triangle and $V\subset H^2(T)$ such that $V\cap P_1(T)=\{0\}$. 
A scaling argument verifies that the best-possible constants $C_0$ and $C_1$ in
\begin{align}\label{eqn:Cs_def}
	\|v\|_{L^2(T)}\leq C_{0}h_T^2\|D^2 v\|_{L^2(T)} \quad\text{and}\quad\|\nabla v\|_{L^2(T)}\leq C_{1}h_T\|D^2
	v\|_{L^2(T)}\qquad\text{for all }v\in V
\end{align}
exclusively depend on the shape of the triangle $T$.
By the Rayleigh-Ritz (or min-max) principle, $\lambda_{\min}=(C_s h_T^s)^{-2}$ is the minimal eigenvalue of the eigenvalue problem
\begin{align}\label{eqn:T_WP}
	(D^2 u, D^2 v)_{L^2(T)}=\lambda b_s(u, v)_{L^2(T)}\qquad\text{ for all }v\in V
\end{align}
with $b_0(\bullet,\bullet) = b(\bullet,\bullet) = (\bullet,\bullet)_{L^2(T)}$ for $s=0$ as in Section~\ref{sec:A
posteriori} and
$b_1(\bullet,\bullet)=(\nabla\bullet,\nabla\bullet)_{L^2(T)}$ for $s=1$ as in plate buckling.
This subsection considers \eqref{eqn:Cs_def}--\eqref{eqn:T_WP} for the four spaces, $V_C\coloneqq H^2_0(T)$ for clamped
boundary conditions, $V_S\coloneqq H^2(T)\cap H^1_0(T)$ for simply-supported boundary conditions, and 
\begin{align*}
	V_V&\coloneqq\{f\in H^2(T)\ :\ f(z)=0 \text{ for all }z\in \V(T)\},\\
	V_M&\coloneqq\{f\in H^2(T)\ :\ f(z)=0 \text{ for all }z\in \V(T)\text{ and }\int_E f\d s=0\text{ for all
	}E\in\E(T)\}.
\end{align*}
The spaces $V_V$ and $V_M$ encode boundary conditions related to the $P_1$ interpolation $I_V:H^2(T)\to P_1(T)$
\cite{BS:MathematicalTheoryFinite2008} and the Morley
interpolation $I_M:H^2(T)\to P_2(T)$ \cite{Mor:TriangularEquilibriumElement1968}
over the triangle $T$.

For the three triangular shapes with $h_T=1$
from Figure~\ref{fig:Triangles}, 
Table~\ref{tab:Triangle_values} displays the principal eigenvalues $\lambda_{\min} = C_s^{-2}$ of \eqref{eqn:T_WP} for
$s\in\{0,1\}$ and the four spaces $V_X$ with $X\in\{C,S,V,M\}$.
In particular, for the Morley-type boundary conditions in $V=V_M$, the constants in the Morley interpolation estimate
read
\begin{align*}
	C_0 = 0.07350005475651561,\qquad
	C_1 = 0.1743250725741249
\end{align*}
for the equilateral triangle and for the right-isosceles triangle
\begin{align*}
	C_0 = 0.045068295511191264,\qquad
	C_1 = 0.16522544473105152.
\end{align*}
Those improve the numerical values in \cite{LSL:OptimalEstimationFujino2019} where it is
suggested that the equilateral
triangle leads to the maximal value $C_0$ amongst all triangles.
\begin{remark}\newcommand{\Const}[1]{C_{#1}}
	Note that the constants $\Const{0},\Const{1}$ for $V=V_{\textrm{M}}$ coincide with the best-possible constants in the local
	Morley interpolation error estimate, i.e., $\Const{0},\Const{1}$ from~\eqref{eqn:Cs_def} for $V=V_{\textrm{M}}$ are
	best-possible in
	\begin{align*}
		\|(1\!-\!I_{\textrm{M}})v\|_{L^2(T)}\leq \Const{0}h_T^2\|D^2 v\|_{L^2(T)},\qquad\|\nabla
		(1\!-\!I_{\textrm{M}})v\|_{L^2(T)}\leq \Const{1}h_T\|D^2
		v\|_{L^2(T)}\qquad\text{for all }v\in H^2(T).
	\end{align*}
	This follows from $V_{\textrm{M}}=(1-I_{\textrm{M}})H^2(T)\subset H^2(T)$, $I_{\textrm{M}}V_{\textrm{M}}=\{0\}$, and
	$\|D^2(1-I_{\textrm{M}})\bullet\|_{L^2(T)}\leq \|D^2 \bullet\|_{L^2(T)}$;
	see~\cite{HSX:ConvergenceOptimalityAdaptive2012,CGH:DiscreteHelmholtzDecomposition2014,Gal:OptimalAdaptiveFEM2015} for details on the Morley
	interpolation $I_{\textrm{M}}$.
\end{remark}

\begin{table}[]
	\centering
	\hspace*{-5em}
	\begin{tabular}{lc|lll}
		$s$
		&	$X$ & equilateral & right-isosceles & $90$-$60$-$30$\\\hline
		\multirow{4}{*}{$0$} 
		& C &	$9804.9449874764568054397360472302$ & $35185.638471713425529039119075$ & $55407.231456202639937488311240$ \\
		& S &	$2770.74747830051377028096946314538$ & $9740.9091034002437236440332688705$ & $15085.180715191686082640833743$ \\
		& V &	$72.942664620393689247350334919683$ & $142.9905816658059570843982031023$ & $120.61629780957915519675157481$ \\
		& M &	$185.10778102243532486304991536227$ & $492.33162470634854919442608899$ & $391.7946452226036199316742042$ \\\hline
		\multirow{4}{*}{$1$} 
		& C &	$146.412905109344790007913155148663$ & $279.14825470003949590687643840$ & $352.12132391711858946305008173$ \\
		& S &	$52.6378901391432459671172853326728$ & $98.696044010893586188344909998761$ & $122.82174365800090725660699910$ \\
		& V &	$9.86394388190996483130098955949938$ & $8.37347027272868454327196355771$ & $8.5380950929242319268917035189$ \\
		& M &	$32.906393793581628348514192294762$ & $36.63077785104390025356999314735$ & $33.45210121650411459840950230$ \\\hline
	\end{tabular}
	\caption{Principal eigenvalues for different boundary conditions in the space $V_X$ with $X\in\{C,S,V,M\}$ on the
	equilateral, right-isosceles, and 90-60-30 triangle of Figure~\ref{fig:Triangles}}
	\label{tab:Triangle_values}
\end{table}

\subsection{Conclusions}%
\label{sub:Conclusions}
The few methods for high-precision reference values in the literature, e.g.,~\cite{BT:HighPrecisionSolutions1999},
highly depend on structural properties of the
domain and do not immediately extend to general domains.
The high-order Argyris FEM demands an adaptive algorithm and the presented hierarchical
Argyris AFEM offers efficient computations of accurate eigenvalues on general polygonal domains. 

The restriction to sufficiently fine initial triangulations and a small bulk parameter $\theta$ for optimal asymptotic convergence rates from 
Theorem~\ref{thmoptimalrates} is \emph{not} visible in the examples.
A coarse initial triangulation or a large bulk parameter $\theta$
close to $1$ increase the preasymptotic regime.
Our overall numerical experience supports our impression that $\theta=0.5$ is a good choice.
After a small number of initial red-refinements to obtain $\T_0$ AFEM enters the asymptotical regime immediately.

While this paper focusses on simple eigenvalues, a generalisation of the theoretical result for multiple eigenvalues or
eigenvalue clusters is possible with the cluster algorithm from~\cite{Gal:AdaptiveFiniteElement2014}.
In all presented applications, however, the AFEM algorithm from Figure~\ref{fig:AFEM} suffices and converges
optimally also towards multiple eigenvalues (e.g., for the square in Subsection~\ref{sub:high_precision}) provided that
the algebraic eigenvalue problem solve is 
sufficiently accurate.

\bibliography{./Bibliography}
\bibliographystyle{alphaabbr}
\end{document}